\documentclass[12pt,a4j, oneside]{amsart}

\usepackage{amsmath,amssymb}
\usepackage{amsthm}
\usepackage[abbrev]{amsrefs}
\usepackage{latexsym}
\usepackage{mathrsfs}
\usepackage{caption}
\usepackage[dvipdfmx]{hyperref}

\theoremstyle{plain}
\newtheorem{thm}{Theorem}[section]
\newtheorem{prop}[thm]{Proposition}
\newtheorem{lem}[thm]{Lemma}
\newtheorem{cor}[thm]{Corollary}
\newtheorem{claim}[thm]{Claim}

\theoremstyle{definition}
\newtheorem{dfn}[thm]{Definition}

\theoremstyle{remark}
\newtheorem{rem}[thm]{Remark}
\newtheorem*{ack}{Acknowledgement}

\makeatletter
\renewcommand{\p@enumii}{}
\makeatother
\numberwithin{equation}{section}

\newcommand{\R}{\mathbb{R}}
\newcommand{\N}{\mathbb{N}}
\newcommand{\C}{\mathbb{C}}
\newcommand{\Hb}{\mathbb{H}}
\newcommand{\Py}{\mathcal{P}}
\newcommand{\PyQ}{\mathcal{Q}}
\newcommand{\X}{\mathcal{X}}

\newcommand{\supp}{\mathop{\mathrm{supp}}}
\newcommand{\diam}{\mathop{\mathrm{diam}}}
\newcommand{\midd}{\mathrel{} \middle| \mathrel{}}

\newcommand{\kf}{d_{\mathrm{KF}}}
\newcommand{\prok}{d_{\mathrm{P}}}
\newcommand{\Lip}{{\mathcal{L}}ip}
\newcommand{\conc}{d_{\mathrm{conc}}}
\newcommand{\haus}{d_{\mathrm{H}}}

\newcommand{\boxarrow}{\xrightarrow{\ \square \ }}
\newcommand{\concarrow}{\xrightarrow{\mathrm{conc}}}
\newcommand{\weakarrow}{\xrightarrow{\mathrm{weak}}}

\title{Convergence of metric transformed spaces}
\author{Daisuke Kazukawa}
\address{Department of Mathematics, Osaka University, Toyonaka, Osaka 560-0032, Japan} 
\email{d-kazukawa@cr.math.sci.osaka-u.ac.jp} 
\thanks{The author is supported by JSPS KAKENHI Grant Number 20J00147}
\date{August 10, 2020}

\begin{document}
\maketitle
\thispagestyle{empty}

\begin{abstract}
We consider the metric transformation of metric measure spaces/pyramids. We clarify the conditions to obtain the convergence of the sequence of transformed spaces from that of the original sequence, and, conversely, to obtain the convergence of the original sequence from that of the transformed sequence, respectively. As an application, we prove that spheres and projective spaces with standard Riemannian distance converge to a Gaussian space and the Hopf quotient of a Gaussian space, respectively, as the dimension diverges to infinity.
\end{abstract}

\setcounter{tocdepth}{1}
\tableofcontents

\section{Introduction}

The geometry and analysis on metric measure spaces have actively been studied. Metric measure spaces typically appear as limit spaces of Riemannian manifolds in the convergence/collapsing theory of Riemannian manifolds. The study of convergence of metric measure spaces is one of central topics in geometric analysis on metric measure spaces. 

Gromov \cite{Grmv}*{Chapter $3\frac{1}{2}_+$} has developed a new convergence theory of metric measure spaces based on the concentration of measure phenomenon studied by L\'evy and V.~Milman \cites{Levy, VMil} (see also \cite{Led}) which is roughly stated as that any 1-Lipschitz function on high-dimensional spaces is close to a constant. Gromov introduced two fundamental concepts of distance functions, the observable distance function $\conc$ and the box distance function $\square$, on the set, say $\X$, of isomorphism classes of metric measure spaces. The box distance function is nearly a metrization of measured Gromov-Hausdorff convergence (precisely the isomorphism classes are little different), while the observable distance function induces a very characteristic topology, called the concentration topology, which admits the convergence of many sequences whose dimensions are unbounded.
The concentration topology is weaker than the box topology and in particular, a measured Gromov-Hausdorff convergence becomes a convergence in the concentration topology. He also introduced a natural compactification, say $\Pi$, of $\X$ with respect to the concentration topology, where the topology on $\Pi$ is called the weak topology. An element $\Py$ of $\Pi$ is called a pyramid and is expressed as a subset of $\X$, so that $\Pi$ is a subset of the power set of $\X$. Under this compactification, we often identify a metric measure space $X$ with a pyramid, say $\Py X$, associated with $X$. We refer to Section 2 for the precise definitions.

The study of the concentration and the weak topologies has been growing in recent years (see \cites{FS, prod, KOS, ellipse, OS, OSu, OY, MML, MMG, ST}). In particular, we have obtained in \cites{ellipse, MML, ST} some nontrivial examples of weak convergent sequences, for example, spheres with the restriction of Euclidean norm, solid ellipses, (projective) Stiefel manifolds with the Frobenius norm, whose dimensions are unbounded. However, in all these examples, the distance function comes from the Euclidean distance. Our final goal in this paper is to give the first nontrivial example of weak convergent sequences of non-Euclidean Riemannian manifolds. For this purpose, we will investigate the convergence of metric transformed spaces.

Let $F \colon [0, +\infty) \to [0, +\infty)$ be a continuous function satisfying the following condition: for any metric space $(X, d_X)$, the function $F \circ d_X$ is a metric on $X$. Such a function $F$ is called a metric preserving function. The concept of metric preserving functions was discovered in 1930s and the study of these functions has been deepened. For example, the function $s \mapsto s/(1 + s)$ is a well-known metric preserving function. More generally, if a function $F \colon [0, +\infty) \to [0, +\infty)$ with $F^{-1}(0) = \{0\}$ is concave, then $F$ is metric-preserving. However, it is also known that metric preserving functions are not necessarily nondecreasing and the class of all metric preserving functions is more complicated. We describe some properties of metric preserving functions in Section 3. 

\begin{dfn}\label{FX}
Let $F \colon [0, +\infty) \to [0, +\infty)$ be a continuous metric preserving function. Given a metric measure space $X = (X, d_X, m_X)$, we define a metric measure space
\begin{equation*}
F(X) := (X, F \circ d_X, m_X).
\end{equation*}
We call $F(X)$ the {\it metric transformed space of $X$ by $F$}. In addition, for a pyramid $\Py$ ($\subset \X$), we define 
\begin{equation*}
F(\Py) := \overline{\bigcup_{X \in \Py} \Py{F(X)}}^{\, \square},
\end{equation*}
where $\overline{\mathcal{Y}}^{\, \square}$ means the $\square$-closure of a family $\mathcal{Y}$ of metric measure spaces.
If $F$ is nondecreasing, $F(\Py)$ is a pyramid and is called the {\it metric transformed pyramid of a pyramid $\Py$ by $F$}.
\end{dfn}

We refer to Proposition \ref{FX:prop} for the proof that $F(\Py)$ is a pyramid if $F$ is nondecreasing. In this proposition, we also show that $F(\Py X) = \Py {F(X)}$ holds for any metric measure space $X$ if $F$ is nondecreasing. However, if not, $F(\Py)$ may not be a pyramid and there exists a metric measure space $X$ such that $F(\Py X) \neq \Py {F(X)}$.

The following theorems are the main results of this paper.

\begin{thm}\label{Py}
Let $F_n, F \colon [0, +\infty) \to [0, +\infty)$, $n = 1, 2, \ldots$, be continuous metric preserving functions. Assume that $F$ is nondecreasing. Then the following {\rm (1)} and {\rm (2)} are equivalent to each other.
\begin{enumerate}
\item For any sequence $\{X_n\}_{n \in \N}$ of metric measure spaces and for any pyramid $\Py$, if $X_n$ converges weakly to $\Py$, then $F_n(X_n)$ converges weakly to $F(\Py)$ as $n \to \infty$.
\item The following three conditions hold.
\begin{enumerate}
\item \label{ptwise} $F_n$ converges pointwise to $F$ as $n \to \infty$.
\item \label{liminf:eq} For any $s \in [0, +\infty)$, 
\begin{equation*}
\lim_{n \to \infty} (F_n(s) - \inf_{s \leq s'} F_n(s')) = 0.
\end{equation*} 
\item \label{limsup:eq}
\begin{equation*}
\limsup_{n \to \infty} \sup F_n \leq \sup F.
\end{equation*}
\end{enumerate}
\end{enumerate}
\end{thm}

\begin{thm}\label{Py_con}
Let $F_n, F \colon [0, +\infty) \to [0, +\infty)$, $n = 1, 2, \ldots$, be continuous metric preserving functions. Assume that $F$ is nondecreasing. Then the following {\rm (1)} and {\rm (2)} are equivalent to each other.
\begin{enumerate}
\item For any sequence $\{X_n\}_{n \in \N}$ of metric measure spaces and for any pyramid $\Py$, if $F_n(X_n)$ converges weakly to $F(\Py)$, then $X_n$ converges weakly to $\Py$ as $n \to \infty$.
\item The following three conditions hold.
\begin{enumerate}
\item $F_n$ converges pointwise to $F$ as $n \to \infty$.
\item For any $s \in [0, +\infty)$, 
\begin{equation*}
\lim_{n \to \infty} (F_n(s) - \inf_{s \leq s'} F_n(s')) = 0.
\end{equation*} 
\setcounter{enumii}{3}
\item \label{incr} $F$ is increasing.
\end{enumerate}
\end{enumerate}
\end{thm}

If all $F_n$ are assumed to be nondecreasing in advance, then we also obtain the following version. We remark that the condition \eqref{liminf:eq} is always true if all $F_n$ are nondecreasing.

\begin{cor}\label{Py:cor}
Let $F_n, F \colon [0, +\infty) \to [0, +\infty)$, $n = 1, 2, \ldots$, be continuous nondecreasing metric preserving functions. Then the following {\rm (A)} and {\rm (B)} hold.
\begin{enumerate}
\item[(A)] The following {\rm (A1)} and {\rm (A2)} are equivalent to each other. 
\begin{enumerate}
\item[(A1)] For any sequence $\{\Py_n\}_{n \in \N}$ of pyramids and for any pyramid $\Py$, if $\Py_n$ converges weakly to $\Py$, then $F_n(\Py_n)$ converges weakly to $F(\Py)$ as $n \to \infty$.
\item[(A2)] \eqref{ptwise} and \eqref{limsup:eq} hold.
\end{enumerate}
\item[(B)] The following {\rm (B1)} and {\rm (B2)} are equivalent to each other. 
\begin{enumerate}
\item[(B1)] For any sequence $\{\Py_n\}_{n \in \N}$ of pyramids and for any pyramid $\Py$, if $F_n(\Py_n)$ converges weakly to $F(\Py)$, then $\Py_n$ converges weakly to $\Py$ as $n \to \infty$.
\item[(B2)] \eqref{ptwise} and \eqref{incr} hold.
\end{enumerate}
\end{enumerate}
\end{cor}

The author obtained in \cite{prod} the similar results in the box and concentration topologies to Theorem \ref{Py}. We investigate some properties in the box and concentration topologies like Theorem \ref{Py_con} and obtain the following result.

\begin{thm}\label{mm_con}
Let $F_n, F \colon [0, +\infty) \to [0, +\infty)$, $n = 1, 2, \ldots$, be continuous metric preserving functions. Then the following {\rm (1)} -- {\rm (3)} are equivalent to each other.
\begin{enumerate}
\item \eqref{ptwise}, \eqref{liminf:eq}, and \eqref{incr} hold.
\item For any sequence $\{X_n\}_{n \in \N}$ of metric measure spaces and for any metric measure space $X$, if $F_n(X_n)$ $\square$-converges to $F(X)$, then $X_n$ $\square$-converges to $X$ as $n \to \infty$.
\item For any sequence $\{X_n\}_{n \in \N}$ of metric measure spaces and for any metric measure space $X$, if $F_n(X_n)$ concentrates to $F(X)$, then $X_n$ concentrates to $X$ as $n \to \infty$.
\end{enumerate}
\end{thm}

\begin{rem}\label{sum_table}
The results obtained in this paper and in \cite{prod} are summarized in the following table (see {\sc Table} 1).

\begingroup
\renewcommand{\arraystretch}{1.2}
\begin{table}[h]
\caption{}
\begin{tabular}{|c|c|}
\hline
`` $X_n \boxarrow X \Rightarrow F_n(X_n) \boxarrow F(X)$ '' &  `` $F_n(X_n) \boxarrow F(X) \Rightarrow X_n \boxarrow X$ '' \\
iff \eqref{ptwise} holds. & iff \eqref{ptwise}, \eqref{liminf:eq}, and \eqref{incr} hold. \\ \hline 
`` $X_n \concarrow X \Rightarrow F_n(X_n) \concarrow F(X)$ '' &  `` $F_n(X_n) \concarrow F(X) \Rightarrow X_n \concarrow X$ '' \\
iff \eqref{ptwise} and \eqref{liminf:eq} hold. & iff \eqref{ptwise}, \eqref{liminf:eq}, and \eqref{incr} hold. \\ \hline 
Assume that $F$ is nondecreasing. & Assume that $F$ is nondecreasing.  \\ 
`` $X_n \weakarrow \Py \Rightarrow F_n(X_n) \weakarrow F(\Py)$ ''  & `` $F_n(X_n) \weakarrow F(\Py) \Rightarrow X_n \weakarrow \Py$ '' \\ 
iff \eqref{ptwise}, \eqref{liminf:eq}, and \eqref{limsup:eq} hold. & iff \eqref{ptwise}, \eqref{liminf:eq}, and  \eqref{incr} hold.  \\ \hline
Assume that all $F_n$ are nondecreasing. & Assume that all $F_n$ are nondecreasing.  \\ 
`` $\Py_n \weakarrow \Py \Rightarrow F_n(\Py_n) \weakarrow F(\Py)$ ''  & `` $F_n(\Py_n) \weakarrow F(\Py) \Rightarrow \Py_n \weakarrow \Py$ '' \\ 
iff \eqref{ptwise} and \eqref{limsup:eq} hold. & iff \eqref{ptwise} and \eqref{incr} hold.  \\ \hline 
\end{tabular}
\end{table}
\endgroup
The top left result is in \cite{prod}*{Corollary 4.4} and the second result on the left-side is in \cite{prod}*{Theorem 1.4} (and see Proposition \ref{Py_ptwise} in this paper).
\end{rem}

\subsection*{Application}
As one of the most important applications of Theorem \ref{Py_con} (or Corollary \ref{Py:cor}), we obtain the weak convergence of spheres and projective spaces with the standard Riemannian distances.

Let $S^n(r)$ be the $n$-dimensional sphere in $\R^{n+1}$ centered at the origin and of radius $r > 0$. We equip $S^n(r)$ with the standard Riemannian distance function and the normalized volume measure. Let $F = \R$, $\C$, or $\Hb$, where $\Hb$ is the algebra of quaternions, and let $d := \dim_\R F$. We consider the Hopf quotient
\begin{equation*}
FP^n(r) := {S^{d(n+1)-1}(r)}/{U^F(1)},
\end{equation*}
where $U^F(1) := \left\{t \in F \midd \|t\| = 1 \right\}$. This is topologically an $n$-dimensional projective space over $F$. We equip $FP^n(r)$ with the quotient metric measure structure of $S^n(r)$. If $F = \C$, then the distance function on ${\C}P^n(r)$ coincides with that induced from the Fubini-Study metric scaled with factor $r$.

\begin{thm}\label{SCP}
Let $\{r_n\}_{n=1}^\infty$ be a given sequence of positive real numbers, and let $\lambda_n := r_n /\sqrt{n}$ {\rm (}resp.~$\lambda_n := r_n /\sqrt{dn}${\rm )}. As $n \to \infty$, $S^n(r_n)$ {\rm (}resp.~$FP^n(r_n)${\rm )} converges weakly to the infinite-dimensional Gaussian space $\Py{\Gamma^\infty_{\lambda^2}}$ {\rm (}resp.~the Hopf quotient $\Py{\Gamma^\infty_{\lambda^2}}/U^F(1)$ of $\Py{\Gamma^\infty_{\lambda^2}}$ {\rm )} if and only if $\lambda_n$ converges to a positive real number $\lambda$.
\end{thm}

We refer to Subsection 6.1 for the definitions of the infinite-dimensional Gaussian space $\Py{\Gamma^\infty_{\lambda^2}}$ and its Hopf quotient $\Py{\Gamma^\infty_{\lambda^2}}/U^F(1)$. We remark that if the distance functions of $S^n(r_n)$ and $FP^n(r_n)$ are induced from the restriction of the Euclidean distance respectively, the weak convergence of $S^n(r_n)$ and $FP^n(r_n)$ has been obtained by Shioya \cites{MML, MMG} and Shioya-Takatsu \cite{ST} (see Theorem \ref{Sthm}).

\begin{ack}
The author would like to thank Professor Takashi Shioya, Professor Takumi Yokota, and Professor Ryunosuke Ozawa for their comments and encouragement.
\end{ack}

\section{Preliminaries}

In this section, we describe the definitions and some properties of metric measure space, the box distance, the observable distance, pyramid, and the weak topology. We use most of these notions along \cite{MMG}. As for more details, we refer to \cite{MMG} and \cite{Grmv}*{Chapter 3$\frac{1}{2}_+$}.

\subsection{Metric measure spaces}
Let $(X, d_X)$ be a complete separable metric space and $m_X$ a Borel probability measure on $X$. We call the triple $(X, d_X, m_X)$ a {\it metric measure space}, or an {\it mm-space} for short. We sometimes say that $X$ is an mm-space, in which case the metric and the measure of $X$ are respectively indicated by $d_X$ and $m_X$.

\begin{dfn}[mm-Isomorphism]
Two mm-spaces $X$ and $Y$ are said to be {\it mm-isomorphic} to each other if there exists an isometry $f \colon \supp{m_X} \to \supp{m_Y}$ such that $f_* m_X = m_Y$, where $f_* m_X$ is the push-forward measure of $m_X$ by $f$. Such an isometry $f$ is called an {\it mm-isomorphism}. Denote by $\mathcal{X}$ the set of mm-isomorphism classes of mm-spaces.
\end{dfn}

Note that an mm-space $X$ is mm-isomorphic to $(\supp{m_X}, d_X , m_X)$. We assume that an mm-space $X$ satisfies 
\begin{equation*}
X = \supp{m_X}
\end{equation*}
unless otherwise stated.

\begin{dfn}[Lipschitz order]
Let $X$ and $Y$ be two mm-spaces. We say that $X$ (\emph{Lipschitz}) \emph{dominates} $Y$ and write $Y \prec X$ if there exists a $1$-Lipschitz map $f \colon X \to Y$ satisfying $f_* m_X = m_Y$. We call the relation $\prec$ on $\X$ the \emph{Lipschitz order}.
\end{dfn}

The Lipschitz order $\prec$ is a partial order relation on $\X$.

\subsection{Box distance and observable distance}
For a subset $A$ of a metric space $(X, d_X)$ and for a real number $r > 0$, we set
\begin{align*}
U_r(A) & :=  \{x \in X \mid d_X(x, A) < r\}, 
\end{align*}
where $d_X(x, A) := \inf_{a \in A} d_X(x, a)$. We sometimes write $U_r^{d_X}(A)$ if we pay attention to the metric $d_X$.

\begin{dfn}[Prokhorov distance]
The \emph{Prokhorov distance} $\prok(\mu, \nu)$ between two Borel probability measures $\mu$ and $\nu$ on a metric space $X$ is defined to be the infimum of $\varepsilon > 0$ satisfying
\begin{equation*}
\mu(U_\varepsilon(A)) \geq \nu(A) - \varepsilon
\end{equation*}
for any Borel subset $A \subset X$. We sometimes write $\prok^{d_X}(\mu, \nu)$ if we pay attention to the metric $d_X$.
\end{dfn}

The Prokhorov metric $\prok$ is a metrization of the weak convergence of Borel probability measures on $X$ provided that $X$ is a separable metric space. Note that if a map $f \colon X \to Y$ between two metric spaces $X$ and $Y$ is $1$-Lipschitz, then we have
\begin{equation}\label{lip_prok}
\prok(f_*\mu, f_* \nu) \leq \prok(\mu, \nu)
\end{equation}
for any two Borel probability measures $\mu$ and $\nu$ on $X$.

\begin{dfn}[Ky Fan metric]
Let $(X, \mu)$ be a measure space and $(Y, d_Y)$ a metric space. For two $\mu$-measurable maps $f,g \colon X \to Y$, we define $\kf^\mu (f, g)$ to be the infimum of $\varepsilon \geq 0$ satisfying
\begin{equation*}
\mu(\{x \in X \mid d_Y(f(x),g(x)) > \varepsilon \}) \leq \varepsilon.
\end{equation*}
The function $\kf^\mu$ is a metric on the set of $\mu$-measurable maps from $X$ to $Y$ by identifying two maps if they are equal to each other $\mu$-almost everywhere. We call $\kf^\mu$ the {\it Ky Fan metric}.
\end{dfn}

\begin{lem}[\cite{MMG}*{Lemma 1.26}]\label{prok_kf}
Let $X$ be a topological space with a Borel probability measure $\mu$ and $Y$ a metric space. For any two Borel measurable maps $f, g \colon X \to Y$, we have
\begin{equation*}
\prok(f_*\mu, g_* \mu) \leq \kf^\mu (f, g).
\end{equation*}
\end{lem}

\begin{dfn}[Parameter]
Let $I := [0,1)$ and let $X$ be an mm-space. A map $\varphi \colon I \to X$ is called a {\it parameter} of $X$ if $\varphi$ is a Borel measurable map such that 
\begin{equation*}
\varphi_\ast \mathcal{L}^1 = m_X,
\end{equation*}
where $\mathcal{L}^1$ is the one-dimensional Lebesgue measure on $I$. 
\end{dfn}

Note that any mm-space has a parameter (see \cite{MMG}*{Lemma 4.2}).

\begin{dfn}[Box distance]
We define the {\it box distance} $\square(X, Y)$ between two mm-spaces $X$ and $Y$ to be the infimum of $\varepsilon \geq 0$ satisfying that there exist parameters $\varphi \colon I \to X$, $\psi \colon I \to Y$, and a Borel subset $I_0 \subset I$ with $\mathcal{L}^1(I_0) \geq 1 - \varepsilon$ such that
\begin{equation*}
|d_X(\varphi(s), \varphi(t)) - d_Y(\psi(s), \psi(t))| \leq \varepsilon
\end{equation*}
for any $s,t \in I_0$.
\end{dfn}

\begin{thm}[\cite{MMG}*{Theorem 4.10}]
The box distance function $\square$ is a complete separable metric on $\mathcal{X}$.
\end{thm}

\begin{lem}[\cite{MMG}*{Proposition 4.12}]\label{mmg4.12}
Let $X$ be a complete separable metric space. For any two Borel probability measures $\mu$ and $\nu$ on $X$, we have
\begin{equation*}
\square((X, \mu), (X, \nu)) \leq 2\prok(\mu, \nu).
\end{equation*}
\end{lem}

The following notion gives one of the conditions that are equivalent to the box convergence.

\begin{dfn}[$\varepsilon$-mm-Isomorphism]
Let $X$ and $Y$ be two mm-spaces and $f \colon X \to Y$ a Borel measurable map. Let $\varepsilon \geq 0$ be a real number. We say that $f$ is an {\it $\varepsilon$-mm-isomorphism} if there exists a Borel subset $X_0 \subset X$ such that 
\begin{enumerate}
\item $m_X(X_0) \geq 1 - \varepsilon$,
\item $|d_X(x, y) - d_Y(f(x),f(y))| \leq \varepsilon$ for any $x, y \in X_0$,
\item $\prok(f_* m_X, m_Y) \leq \varepsilon$.
\end{enumerate}
We call $X_0$ a {\it nonexceptional domain} of $f$.
\end{dfn}
It is easy to see that, for a 0-mm-isomorphism $f \colon X \to Y$, there is an mm-isomorphism $\hat{f} \colon X \to Y$ that is equal to $f$ $m_X$-a.e.~on $X$.

\begin{lem}[\cite{MMG}*{Lemma 4.22}]\label{lem4.22}
Let $X$ and $Y$ be two mm-spaces.
\begin{enumerate}
\item If there exists an $\varepsilon$-mm-isomorphism $f \colon X \to Y$, then $\square(X, Y) \leq 3\varepsilon$.
\item If $\square(X, Y) < \varepsilon$, then there exists a $3\varepsilon$-mm-isomorphism $f\colon X \to Y$.
\end{enumerate}
\end{lem}

Given an mm-space $X$ and a parameter $\varphi \colon I \to X$ of $X$, we set
\begin{equation}\label{pull_lip1}
\varphi^* \Lip_1(X) := \{ f \circ \varphi \mid f \colon X \to \R \text{ is $1$-Lipschitz} \}.
\end{equation}
Note that $\varphi^* \Lip_1(X)$ consists of Borel measurable functions on $I$.

\begin{dfn}[Observable distance] 
We define the {\it observable distance} $\conc(X, Y)$ between two mm-spaces $X$ and $Y$ by
\begin{equation*}
\conc(X, Y) := \inf_{\varphi, \psi} \haus(\varphi^* \Lip_1(X), \psi^* \Lip_1(Y)),
\end{equation*}
where $\varphi \colon I \to X$ and $\psi \colon I \to Y$ run over all parameters of $X$ and $Y$ respectively, and $\haus$ is the Hausdorff distance with respect to the metric $\kf^{\mathcal{L}^1}$. We say that a sequence $\{X_n\}_{n \in \N}$ of mm-spaces  {\it concentrates} to an mm-space $X$ if $X_n$ $\conc$-converges to $X$ as $n \to \infty$.
\end{dfn}

\begin{prop}[\cite{MMG}*{Proposition 5.5}]
For any two mm-spaces $X$ and $Y$, we have
\begin{equation*}
\conc(X, Y) \leq \square(X, Y).
\end{equation*}
\end{prop}

\begin{thm}[\cite{MMG}*{Theorem 5.13}]
The observable distance function $\conc$ is a metric on $\mathcal{X}$.
\end{thm}

The basic lemmas used in this paper are listed as follows.

\begin{lem}[\cite{MMG}*{Corollary 4.48}]\label{fin_dim_app}
For any mm-space $X$, there exist $1$-Lipschitz maps $\Phi_N \colon X \to (\R^N, \|\cdot\|_\infty)$, $N = 1, 2, \ldots$, such that
\begin{equation*}
\lim_{N \to \infty} \square(X, (\R^N, \|\cdot\|_\infty, {\Phi_N}_* m_X)) = 0,
\end{equation*}
where $\|x - y\|_\infty := \max_{1 \leq i \leq N} |x_i - y_i|$ for any $x, y \in \R^N$.
\end{lem}

\begin{dfn}[$1$-Lipschitz up to an additive error]
Let $X$ be an mm-space and $Y$ be a metric space. A map $f \colon X \to Y$ is said to be 1-{\it Lipschitz up to} ({\it an additive error}) {\it $\varepsilon \geq 0$} if there exists a Borel subset $X_0 \subset X$ such that
\begin{enumerate}
\item $m_X(X_0) \geq 1 - \varepsilon$,
\item $d_Y(f(x), f(x')) \leq d_X(x, x') + \varepsilon$ for any $x, x' \in X_0$.
\end{enumerate}
We call $X_0$ a {\it nonexceptional domain} of $f$.
\end{dfn}

\begin{lem}[\cite{MMG}*{Lemma 5.4}]\label{MWext}
If a Borel measurable function $f \colon X \to (\R^N, \|\cdot\|_\infty)$ on an mm-space $X$ is $1$-Lipschitz up to an additive error $\varepsilon \geq 0$, then there exists a $1$-Lipschitz function $\tilde{f} \colon X \to (\R^N, \|\cdot\|_\infty)$ such that
\begin{equation*}
\kf^{m_X}(f, \tilde{f}) \leq \varepsilon.
\end{equation*}
\end{lem}

\begin{thm}[\cite{MMG}*{Theorem 4.35}]\label{thm4.35}
Let $X$, $Y$, $X_n$, and $Y_n$ be mm-spaces, $n = 1,2,\ldots$. If $X_n$ and $Y_n$ $\square$-converge to $X$ and $Y$ respectively as $n \to \infty$ and if $X_n \prec Y_n$ for any $n$, then $X \prec Y$.
\end{thm}

\begin{lem}[\cite{MMG}*{Lemma 6.10}]\label{lem6.10}
\begin{enumerate}
\item If a sequence $\{X_n\}_{n \in \N}$ of mm-spaces $\square$-converges to an mm-space $X$ and if $X$ dominates an mm-space $Y$, then there exists a sequence $\{Y_n\}_{n \in \N}$ of mm-spaces $\square$-converging to $Y$ such that $X_n$ dominates $Y_n$ for each $n$.
\item If two sequences $\{X_n\}_{n \in \N}$ and $\{Y_n\}_{n \in \N}$ of mm-spaces $\square$-converge and if $X_n$ and $Y_n$ are both dominated by an mm-space $\tilde{Z}_n$ for each $n$, then there exists a sequence of mm-spaces $Z_n$ such that $X_n,Y_n \prec Z_n \prec \tilde{Z}_n$ and $\{Z_n\}_{n \in \N}$ has a $\square$-convergent subsequence.
\end{enumerate}
\end{lem}

\subsection{Pyramid}

\begin{dfn}[Pyramid] \label{Py:dfn}
A subset $\Py \subset \X$ is called a \emph{pyramid} if it satisfies the following {\rm(1) -- (3)}.
\begin{enumerate}
\item If $X \in \Py$ and if $Y \prec X$, then $Y \in \Py$.
\item For any $X, X' \in \Py$, there exists $Y \in \Py$ such that $X \prec Y$ and $X' \prec Y$.
\item $\Py$ is nonempty and $\square$-closed.
\end{enumerate}
We denote the set of pyramids by $\Pi$. Note that Gromov's definition of a pyramid is only by (1) and (2). The condition (3) is added in \cite{MMG} for the Hausdorff property of $\Pi$.

For an mm-space $X$, we define
\begin{equation*}
\Py X := \left\{X' \in \X \midd X' \prec X \right\},
\end{equation*}
which is a pyramid. We call $\Py X$ the \emph{pyramid associated with $X$}.
\end{dfn}

We observe that $X \prec Y$ if and only if $\Py X \subset \Py Y$. 

\begin{dfn}[Weak convergence]
Let $\Py$, $\Py_n \in \Pi$, $n = 1, 2, \ldots$. We say that $\Py_n$ \emph{converges weakly to} $\Py$ as $n \to \infty$ if the following (1) and (2) are both satisfied.
\begin{enumerate}
\item For any mm-space $X \in \Py$, we have
\begin{equation*}
\lim_{n \to \infty} \square(X, \Py_n) = 0.
\end{equation*}
\item For any mm-space $X \in \X \setminus \Py$, we have
\begin{equation*}
\liminf_{n \to \infty} \square(X, \Py_n) > 0.
\end{equation*}
\end{enumerate}
\end{dfn}

\begin{thm}\label{Py:thm}
There exists a metric, denoted by $\rho$, on $\Pi$ such that the following {\rm (1) -- (4)} hold.
\begin{enumerate}
\item $\rho$ is compatible with weak convergence.
\item The map
$\iota \colon \X \ni X \mapsto \Py X \in \Pi$
is a $1$-Lipschitz topological embedding map with respect to
$\conc$ and $\rho$.
\item $\Pi$ is $\rho$-compact.
\item $\iota(\X)$ is $\rho$-dense in $\Pi$.
\end{enumerate}
\end{thm}
In particular, $(\Pi,\rho)$ is a compactification of $(\X,\conc)$. We often identify $X$ with $\Py X$, and we say that a sequence of mm-spaces \emph{converges weakly} to a pyramid if the associated pyramid converges weakly.

\begin{lem}[\cite{MMG}*{Lemma 7.14}]\label{approximation}
For any pyramid $\Py$, there exists a sequence $\{Y_m\}_{m \in \N}$ of mm-spaces such that
\begin{equation*}
Y_1 \prec Y_2 \prec \cdots \prec Y_m \prec \cdots \quad \text{ and } \quad \overline{\bigcup_{m = 1}^\infty \Py{Y_m}}^{\, \square} = \Py.
\end{equation*}
\end{lem}

Such a sequence $\{Y_m\}_{m \in \N}$ is called an \emph{approximation of $\Py$}. We see that $Y_m$ converges weakly to $\Py$ as $m \to \infty$ and that $Y_m \in \Py$ for all $m$. 

\begin{lem}[cf.~\cite{Grmv}*{3$\frac{1}{2}$.15.}]\label{chara_ass_py}
Let $\Py$ be a pyramid. The following {\rm (1)} and {\rm (2)} are equivalent to each other.
\begin{enumerate}
\item $\Py \in \iota(\X)$, i.e., there exists an mm-space $X$ such that $\Py = \Py X$.
\item $\Py$ is $\square$-compact.
\end{enumerate}
\end{lem}

\begin{proof}
We first prove `(1) $\Rightarrow$ (2)'. Take any mm-space $X$ and prove that $\Py X$ is $\square$-compact. We take any real number $\varepsilon > 0$. By \cite{MMG}*{Lemma 4.28}, it is sufficient to prove that there exists a real number $\Delta(\varepsilon) > 0$ such that for any $Y \in \Py X$ we have a finite net $\mathcal{N} \subset Y$ such that
\begin{equation*}
m_Y(U_\varepsilon(\mathcal{N})) \geq 1 -\varepsilon, \quad \# \mathcal{N} \leq \Delta(\varepsilon), \quad \text{and} \quad \diam{\mathcal{N}} \leq \Delta(\varepsilon).
\end{equation*}
We find a finite net $\mathcal{N} \subset X$ with $m_X(U_\varepsilon(\mathcal{N})) \geq 1 -\varepsilon$. Note that the existence of such $\mathcal{N}$ follows from the separability of $X$. We define
\begin{equation*}
\Delta(\varepsilon) := \max\{\#\mathcal{N}, \diam{\mathcal{N}}\}.
\end{equation*}
Take any mm-space $Y \in \Py X$ and fix it. There exists a $1$-Lipschitz map $f \colon X \to Y$ such that $f_* m_X = m_Y$. Then the finite net $f(\mathcal{N})$ of $Y$ satisfies
\begin{align*}
&\# f(\mathcal{N}) \leq \# \mathcal{N} \leq \Delta(\varepsilon), \\
&\diam{f(\mathcal{N})} \leq \diam{\mathcal{N}} \leq \Delta(\varepsilon), \\
&m_Y(U_\varepsilon(f(\mathcal{N}))) \geq m_X(U_\varepsilon(\mathcal{N})) \geq 1-\varepsilon.
\end{align*}
Thus $\Py X$ is $\square$-compact.

We next prove `(2) $\Rightarrow$ (1)'. Let $\{Y_m\}_{m \in \N}$ be an approximation of $\Py$. Note that $Y_m$ converges weakly to $\Py$. Since $\Py$ is $\square$-compact, the sequence $\{Y_m\}_{m \in \N}$ has a $\square$-convergent subsequence. Let $X$ be a limit mm-space of a $\square$-convergent subsequence of $\{Y_m\}_{m \in \N}$. Since a $\square$-convergence becomes a weak convergence, we have $\Py = \Py X$. This completes the proof.
\end{proof}

\section{Metric preserving functions}
In this section, we recall some properties of metric preserving function. The notion of metric preserving functions was discovered in 1930s, and the study of these functions has been deepened. We refer to \cite{C} for a survey of results on metric preserving functions.

\subsection{Metric preserving functions}

\begin{dfn}[Metric preserving function]
A function $F \colon [0, +\infty) \to [0, +\infty)$ is called a {\it metric preserving function} provided that for any metric space $(X, d_X)$, the function $F \circ d_X$ is a metric on $X$.
\end{dfn}

\begin{lem}[\cite{C}*{Propositions 2.1, 2.3}]\label{Kelly}
Let $F \colon [0, +\infty) \to [0, +\infty)$ be a function. Then the following {\rm (1)} and {\rm (2)} hold.
\begin{enumerate}
\item If $F$ is a metric preserving function, then $F$ is subadditive {\rm (}i.e., 
\begin{equation*}
F(s + t) \leq F(s)  + F(t)
\end{equation*}
for any $s, t \geq 0${\rm )} and $F^{-1}(0) = \{0\}$.
\item If $F$ is subadditive and nondecreasing and fulfills $F^{-1}(0) = \{0\}$, then $F$ is a metric preserving function. In particular, if $F$ is a concave function with $F^{-1}(0) = \{0\}$, then $F$ is a metric preserving function.
\end{enumerate}
\end{lem}

\begin{rem}
There are many examples of metric preserving functions that are not nondecreasing. For example,
\begin{align*}
F(s) & := \left\{ \begin{array}{ll} s & \text{ if } s \in [0, 2), \\ 4 - s & \text{ if } s \in [2, 3), \\ 1 & \text{ if } s \in [3, +\infty), \end{array} \right. \\
G(s) & := \left\{ \begin{array}{ll} s & \text{ if } s \in [0, 1), \\ \displaystyle \frac{1+s+\sin^2(s-1)}{2s} & \text{ if } s \in [1, +\infty). \end{array} \right.
\end{align*}
\end{rem}

\begin{prop}[cf.~\cite{C}*{Propositions 2.6}]\label{cor_mpf}
Let $F \colon [0, +\infty) \to [0, +\infty)$ be a metric preserving function. Then, for any $s, t  \geq 0$, we have
\begin{enumerate}
\item $|F(s) - F(t)| \leq F(|s-t|)$,
\item $F(s) \leq 2 F(t)$ if $s \leq 2t$.
\end{enumerate}
\end{prop}

\begin{thm}[\cite{C}*{Theorem 3.4}]\label{Fconti}
Let $F \colon [0, +\infty) \to [0, +\infty)$ be a metric preserving function. Then the following conditions are equivalent to each other.
\begin{enumerate}
\item $F$ is continuous.
\item $F$ is continuous at $0$.
\item $F$ is uniformly continuous.
\item For any metric space $(X, d_X)$, the topologies induced by $d_X$ and $F \circ d_X$ coincide with each other. 
\end{enumerate}
\end{thm}

Note that if $F$ is discontinuous, then $F \circ d_X$ gives the discrete topology on $X$.

\begin{prop}[cf.~\cite{prod}*{Proposition 3.10}]\label{Fcplt}
Let $F \colon [0, +\infty) \to [0, +\infty)$ be a metric preserving function. If a metric space $(X, d_X)$ is complete, then so is $(X, F \circ d_X)$.
\end{prop}

\subsection{Convergence of metric preserving functions}

In this subsection, we describe some properties for a sequence of metric preserving functions. In particular, we show some conditions that are equivalent to \eqref{ptwise}, \eqref{liminf:eq}, and \eqref{limsup:eq} respectively.

\begin{lem}\label{ptcpt}
Let $F_n, F \colon [0, +\infty) \to [0, +\infty)$, $n = 1, 2, \ldots$, be metric preserving functions. If $F$ is continuous and if $F_n$ converges pointwise to $F$, then $F_n$ converges uniformly to $F$ on compact sets.
\end{lem}

\begin{proof}
We take any compact set $K \subset [0, +\infty)$ and any real number $\varepsilon > 0$. Let us prove that
\begin{equation}\label{ptcpt_eq}
\sup_{s \in K}|F_n(s) - F(s)| \leq 7\varepsilon
\end{equation}
holds for every sufficiently large $n$. By the continuity of $F$, there exists a real number $\delta > 0$ such that $F(\delta) \leq \varepsilon$. By the compactness of $K$, we find finite points $\{s_i\}_{i=1}^k$ in $K$ such that
\begin{equation*}
K \subset \bigcup^k_{i=1} U_\delta(s_i).
\end{equation*}
Let $N \in \N$ be a number such that
\begin{equation*}
\max_{i=1, \ldots, k} |F_n(s_i) - F(s_i)| \leq \varepsilon \quad \text{and} \quad |F_n(\delta) - F(\delta)| \leq \varepsilon
\end{equation*}
hold for all $n \geq N$. Given a fixed point $s \in K$, we find $i \in \{1, \ldots, k\}$ such that $s \in U_\delta(s_i)$. By Proposition \ref{cor_mpf}, we have
\begin{align*}
&|F_n(s) - F(s)| \\
&\leq |F_n(s_i) - F(s_i)| + F_n(|s - s_i|) + F(|s - s_i|) \\
&\leq |F_n(s_i) - F(s_i)| + 2F_n(\delta) + 2F(\delta) \\
&\leq |F_n(s_i) - F(s_i)| + 2|F_n(\delta) - F(\delta)| + 4F(\delta) \\
&\leq 7\varepsilon
\end{align*}
for every $n \geq N$. Thus we obtain \eqref{ptcpt_eq}. This completes the proof.
\end{proof}

Given a function $F \colon [0, +\infty) \to [0, +\infty)$, we set
\begin{equation*}
I_F(s) : = F(s) - \inf_{s \leq s'} F(s'), \quad s \in [0, +\infty).
\end{equation*}
Note that $I_F \geq 0$ and that $I_F \equiv 0$ if and only if $F$ is nondecreasing.

\begin{lem}\label{liminf:lem}
Let $F_n \colon [0, +\infty) \to [0, +\infty)$ be functions, $n = 1, 2, \ldots$. Assume that $F_n$  converges uniformly to a continuous function $F$ on compact sets. Then the following {\rm (1)} -- {\rm (3)} are equivalent to each other. 
\begin{enumerate}
\item For any $s > 0$,
\begin{equation*}
\lim_{n \to \infty} I_{F_n}(s) = 0
\end{equation*}
holds {\rm (}i.e., \eqref{liminf:eq} holds{\rm )}.
\item For any $D > 0$,
\begin{equation*}
\lim_{n \to \infty} \sup_{s \in [0, D]} I_{F_n}(s) = 0
\end{equation*}
holds.
\item $F$ is nondecreasing and, for any sequence $s_n \to \infty$,
\begin{equation*}
\liminf_{n \to \infty} F_n(s_n) \geq \sup{F} \ \left( = \lim_{n \to \infty} F(s_n) \right)
\end{equation*}
holds.
\end{enumerate}
\end{lem}

\begin{proof}
`$(2) \Rightarrow (1)$' is obvious. We verify `$(1) \Rightarrow (3)$' and `$(3) \Rightarrow (2)$'.

Assume (1). Take any two real numbers $s, s'$ with $0 \leq s \leq s'$. Then
\begin{equation*}
F(s) = \lim_{n \to \infty} F_n(s) \leq \lim_{n \to \infty} (F_n(s') + I_{F_n}(s)) = F(s'),
\end{equation*}
which implies that $F$ is nondecreasing. We take any sequence $s_n \to \infty$ and any real number $s \geq 0$. For every sufficiently large $n$, since $s \leq s_n$, we see that
\begin{equation*}
F_n(s) \leq F_n(s_n) + I_{F_n}(s).
\end{equation*}
Thus we have
\begin{equation*}
F(s) = \lim_{n \to \infty} F_n(s) \leq \liminf_{n \to \infty} (F_n(s_n) + I_{F_n}(s)) =\liminf_{n \to \infty} F_n(s_n),
\end{equation*}
which implies
\begin{equation*}
\sup{F} \leq \liminf_{n \to \infty} F_n(s_n).
\end{equation*}
The proof of `$(1) \Rightarrow (3)$' is completed.

We next prove `$(3) \Rightarrow (2)$'. Suppose that (2) does not hold. Then there exists $D > 0$ such that
\begin{equation*}
\eta := \limsup_{n \to \infty} \sup_{s \in [0, D]} I_{F_n}(s) > 0.
\end{equation*}
Taking a subsequence of $n$, we can assume that $\sup_{s \in [0, D]} I_{F_n}(s) \to \eta$. Thus, for every sufficiently large $n$, we have
\begin{equation*}
\sup_{s \in [0, D]} I_{F_n}(s) > \frac{\eta}{2}.
\end{equation*}
Then, there exist two real numbers $s_n, s'_n$ with $0 < s_n \leq \min\{s'_n, D\}$ such that
\begin{equation*}
F_n(s_n) > F_n(s'_n) + \frac{\eta}{2}.
\end{equation*}
Taking a subsequence again, we are able to assume that at least one of the following two situations occurs.
\begin{itemize}
\item $s_n$ and $s'_n$ converge to real numbers $s_\infty$ and $s'_\infty$ respectively.
\item $s_n$ converges to a real number $s_\infty$ and $s'_n$ diverges to $+\infty$.
\end{itemize}
If the first situation occurs, then we have
\begin{equation*}
s_\infty \leq s'_\infty \quad \text{and} \quad F(s_\infty) \geq F(s'_\infty) + \frac{\eta}{2}.
\end{equation*}
In fact, since $F_n$ converges uniformly to $F$ on compact sets, we have $F_n(s_n) \to F(s)$ if $s_n \to s$. However this contradicts the monotonicity of $F$. If the second situation occurs, then we have
\begin{equation*}
F(s_\infty) \geq \liminf_{n \to \infty} F_n(s'_n) + \frac{\eta}{2},
\end{equation*}
which contradicts $\liminf_{n \to \infty} F_n(s'_n) \geq \sup{F}$. Therefore we obtain `$(3) \Rightarrow (2)$'. 

The proof of this lemma is completed.
\end{proof}

\begin{rem}\label{compcond}
Under the setting of Lemma \ref{liminf:lem}, we consider the following other conditions.
\begin{enumerate}
\item[(i)] The functions $F_n$ are nondecreasing for all $n \in \N$.
\item[(ii)] $\lim_{n\to\infty} \sup_{s \geq 0} I_{F_n}(s) = 0$.
\item[(iii)] The function $F$ is nondecreasing.
\end{enumerate}
It is easy to see that `(i) $\Rightarrow$ (ii) $\Rightarrow$ \eqref{liminf:eq} $\Rightarrow$ (iii)'. On the other hand, the converse of each of these implications does not hold, even in the class of metric preserving functions. In fact, we show the following examples. 
We define functions $F_n^1$, $F_n^2$, and $F_n^3$, $n = 1, 2, \ldots$, by
\begin{align*}
F_n^1(s) & := \left\{ \begin{array}{ll} s & \text{ if } s \in [0, 2), \\  4 - s & \text{ if } s \in [2, 2+n^{-1}), \\ 2 - n^{-1} & \text{ if } s \in [2 + n^{-1}, +\infty). \end{array} \right.
\end{align*}
\begin{align*}
F_n^2(s) & := \left\{ \begin{array}{ll} s & \text{ if } s \in [0, 2), \\ 2 & \text{ if } s \in [2, n+2), \\  s - n & \text{ if } s \in [n+2, n+3), \\ n + 6 - s & \text{ if } s \in [n+3, n+4), \\ 2 & \text{ if } s \in [n+4, +\infty). \end{array} \right.
\end{align*}
\begin{align*}
F_n^3(s) & := \left\{ \begin{array}{ll} s & \text{ if } s \in [0, 2), \\ 2 & \text{ if } s \in [2, n+2), \\  n + 4 - s & \text{ if } s \in [n+2, n+3), \\ 1 & \text{ if } s \in [n+3, +\infty). \end{array} \right.
\end{align*}
These functions are continuous metric preserving functions and converge to the concave function $\min\{s, 2\}$ as $n \to \infty$. $\{F_n^1\}$, $\{F_n^2\}$, and $\{F_n^3\}$ are counterexamples of `(ii) $\Rightarrow$ (i)', `\eqref{liminf:eq} $\Rightarrow$ (ii)', and `(iii) $\Rightarrow$ \eqref{liminf:eq}' respectively.
\end{rem}

\begin{lem}\label{limsup:lem}
Let $F_n \colon [0, +\infty) \to [0, +\infty)$ be functions, $n = 1, 2, \ldots$. Assume that $F_n$  converges uniformly to a continuous function $F$ on compact sets. Then the following {\rm (1)} -- {\rm (3)} are equivalent to each other. 
\begin{enumerate}
\item 
\begin{equation*}
\limsup_{n \to \infty} \sup{F_n} \leq \sup{F}
\end{equation*}
holds {\rm (}i.e., \eqref{limsup:eq} holds{\rm )}.
\item 
\begin{equation*}
\limsup_{n \to \infty} \sup{F_n} = \sup{F}
\end{equation*}
holds.
\item For any sequence $s_n \to \infty$,
\begin{equation*}
\limsup_{n \to \infty} F_n(s_n) \leq \sup{F}
\end{equation*}
holds. 
\end{enumerate}
\end{lem}

\begin{proof}
We first verify that $\liminf_{n \to \infty} \sup{F_n} \geq \sup{F}$ is always true. For any $s \geq 0$, we have
\begin{equation*}
F(s) = \lim_{n \to \infty} F_n(s) \leq \liminf_{n \to \infty} \sup{F_n},
\end{equation*}
which implies
\begin{equation*}
\sup{F} \leq \liminf_{n \to \infty} \sup{F_n}.
\end{equation*}
Therefore we obtain `(1) $\Leftrightarrow$ (2)'. Moreover, `(1) $\Rightarrow$ (3)' is trivial.

We verify `(3) $\Rightarrow$ (1)'. We first assume that $\sup{F_n} < +\infty$ for all $n$. We take any real number $\varepsilon > 0$. There exists a sequence $\{s_n\}$ of positive real numbers such that
\begin{equation*}
\sup{F_n} \leq F_n(s_n) + \varepsilon
\end{equation*}
for every $n$. Taking a subsequence, we can assume that $\{s_n\}$ converges to a real number $s_\infty$ or it diverges to infinity. If $s_n \to s_\infty$, then we have
\begin{equation*}
\limsup_{n \to \infty} \sup{F_n} \leq  F(s_\infty) + \varepsilon \leq \sup{F} + \varepsilon.
\end{equation*}
If $s_n \to +\infty$, then we have
\begin{equation*}
\limsup_{n \to \infty} \sup{F_n} \leq  \limsup_{n \to \infty} F_n(s_n) + \varepsilon \leq \sup{F} + \varepsilon.
\end{equation*}
Thus, as $\varepsilon \to 0$, we obtain (1). In the case that $\sup{F_{n_i}} = +\infty$ for some subsequence $\{n_i\}$, for any real number $M > 0$, there exists a sequence $\{s_i\}$ of positive real numbers such that
\begin{equation*}
F_{n_i}(s_i) > M
\end{equation*}
for every $i$. In the same discussion as above, taking a subsequence of $\{n_i\}$, we obtain
\begin{equation*}
M < \limsup_{i \to \infty} F_{n_i}(s_i) \leq \sup{F},
\end{equation*}
which implies $\sup{F} = +\infty$. Thus we obtain (1) in general.

The proof of this lemma is completed.
\end{proof}

\begin{prop}\label{ptwise_con}
Let $F_n, F \colon [0,+\infty) \to [0, +\infty)$, $n = 1, 2, \ldots$, be continuous metric preserving  functions. Assume that, for any sequence $\{s_n\} \subset [0, +\infty)$ and any $s > 0$, if $F_n(s_n)$ converges to $F(s)$, then $s_n$ converges to $s$. Then $F_n$ converges pointwise to $F$.
\end{prop}

\begin{proof}
We take any $s > 0$ and fix it. We first prove 
\begin{equation}\label{conv_mpf_liminf:eq}
F(s) \leq \liminf_{n \to \infty} F_n(s).
\end{equation}
We set $\alpha := \liminf_{n \to \infty} F_n(s)$ and suppose that $\alpha < F(s)$. There exists a subsequence $\{n_i\}_i$ of $n$ such that $F_{n_i}(s) \to \alpha$ as $i \to \infty$. For each $n$, if $\alpha < F_n(s)$, then there exists a real number $s'_n > 0$ such that
\begin{equation*}
F_n(s'_n)  = \alpha
\end{equation*}
by the intermediate value theorem. We set a sequence
\begin{equation*}
s_n := \left\{ \begin{array}{ll} s'_n & \text{ if } n \neq n_i \text{ and if } \alpha < F_n(s), \\ s & \text{ otherwise. } \end{array} \right. 
\end{equation*}
Taking the definition of $\alpha$ into account, we see that $F_n(s_n)$ converges to $\alpha$ as $n \to \infty$. By $\alpha < F(s)$ and by the intermediate value theorem, there exists a real number $\beta > 0$ such that $\alpha = F(\beta)$. Thus, by the assumption of this proposition, we have $s_n \to \beta$ as $n \to \infty$. Since $\{s_n\}$ has a subsequence consisting only of $s$, we have $\beta = s$, which contradicts $\alpha < F(s)$. Thus we have \eqref{conv_mpf_liminf:eq}.

We next prove 
\begin{equation}\label{conv_mpf_limsup:eq}
\limsup_{n \to \infty} F_n(s) \leq F(s).
\end{equation}
Suppose that $F_n(s) > F(s) + \eta$ for every sufficiently large $n$ and for some real number $\eta > 0$. For every sufficiently large $n$, there exists a real number $s_n > 0$ such that $s_n \leq s$ and $F_n(s_n) = F(s)$. By the assumption of this proposition, $s_n$ converges to $s$. Since $F_n$ is a metric preserving function, for every sufficiently large $n$, we have
\begin{equation*}
0 < \eta < F_n(s) - F(s) =  F_n(s) - F_n(s_n) \leq F_n(s - s_n).
\end{equation*}
Let $\eta' := \min\{F(s), \eta\} > 0$. For every sufficiently large $n$, there exists a real number $t_n$ such that $t_n \leq s - s_n$ and $F_n(t_n) = \eta'$. We see that $t_n$ converges to $0$. On the other hand, since a real number $t > 0$ such that $F(t) = \eta'$ is also found, it is follows from the assumption that $t_n$ converges to $t$. This is a contradiction. Thus we obtain \eqref{conv_mpf_limsup:eq}. The proof is completed.
\end{proof}

\section{Weak convergence of metric transformed pyramids}

The goal in this section is to prove Theorem \ref{Py} and Theorem \ref{Py_con}.

We review the definitions of $F(X)$ and $F(\Py)$ in Definition \ref{FX}. Let $F \colon [0, +\infty) \to [0, +\infty)$ be a continuous metric preserving function. Given an mm-space $X$ and a pyramid $\Py$, we define 
\begin{equation*}
F(X) := (X, F \circ d_X, m_X) \quad \text{and} \quad F(\Py) := \overline{\bigcup_{X \in \Py} \Py{F(X)}}^{\, \square}.
\end{equation*}
By Theorem \ref{Fconti} and Proposition \ref{Fcplt}, $F(X)$ is an mm-space for a given mm-space $X$.

\begin{prop}\label{FX:prop}
Let $F \colon [0, +\infty) \to [0, +\infty)$ be a continuous metric preserving function. If $F$ is nondecreasing, then $F(\Py)$ is a pyramid for any pyramid $\Py$ and $F(\Py X) = \Py F(X)$ holds for any mm-space $X$.
\end{prop}

\begin{proof}
Let $\Py$ be a pyramid. We verify that $F(\Py)$ is a pyramid. It is obvious that $F(\Py)$ is nonempty and $\square$-closed.

We check the condition (1) of Definition \ref{Py:dfn}. Assume $Y \in F(\Py)$ and $Y' \prec Y$. By the definition of $F(\Py)$, there exist mm-spaces $X_n \in \Py$ and $Y_n \in \X$, $n = 1, 2, \ldots$, such that $F(X_n)$ dominates $Y_n$ for every $n$ and $\square(Y_n, Y) \to 0$ as  $n \to \infty$. By Lemma \ref{lem6.10} (1), there exists a sequence $\{Y'_n\}_{n \in \N}$ of mm-spaces such that $Y_n$ dominates $Y'_n$ for every $n$ and $\square(Y'_n, Y') \to 0$ as  $n \to \infty$. Since $F(X_n)$ dominates $Y'_n$ (i.e., $Y'_n \in \Py F(X_n)$) for every $n$, we have $Y' \in F(\Py)$. 

We next check the condition (2) of Definition \ref{Py:dfn}. Take any two mm-spaces $Y, Y' \in F(\Py)$. By the definition of $F(\Py)$, there exist mm-spaces $X_n, X'_n \in \Py$ and $Y_n, Y'_n \in \X$, $n = 1, 2, \ldots$, such that $F(X_n)$ (resp.~$F(X'_n)$) dominates $Y_n$ (resp.~$Y'_n$) for every $n$ and $\square(Y_n, Y) \to 0$, $\square(Y'_n, Y') \to 0$ as $n \to \infty$. By $X_n, X'_n \in \Py$, there exists $\widetilde{X}_n \in \Py$ such that $\widetilde{X}_n$ dominates both $X_n$ and $X'_n$. Since $F$ is nondecreasing, we see that $F(\widetilde{X}_n)$ dominates both $F(X_n)$ and $F(X'_n)$, which implies that $F(\widetilde{X}_n)$ dominates both $Y_n$ and $Y'_n$. By Lemma \ref{lem6.10} (2), there exists a sequence $\{Z_n\}_{n \in \N}$ of mm-spaces such that $Y_n, Y'_n \prec Z_n \prec F(\widetilde{X}_n)$ for every $n$ and $\{Z_n\}_{n \in \N}$ has a $\square$-convergent subsequence. Let $Z$ be a limit  space of $\square$-convergent subsequence of $\{Z_n\}_{n \in \N}$. Since $\{Z_n\}_{n \in \N} \subset F(\Py)$, we see that $Z \in F(\Py)$ and, by Theorem \ref{thm4.35}, we have $Y, Y' \prec Z$. Thus $F(\Py)$ satisfies the condition (2), so that $F(\Py)$ is a pyramid.

We prove that $F(\Py X) = \Py F(X)$ for any mm-space $X$. We take an mm-space $X$ and fix it. Since $X \in \Py X$, we have $F(\Py X) \supset \Py F(X)$. Since $F$ is nondecreasing, if $Y \in \Py X$, then $\Py F(Y) \subset \Py F(X)$, which leads to $F(\Py X) \subset \Py F(X)$. This completes the proof.
\end{proof}

\subsection{Proof of Theorem \ref{Py}}

\begin{prop}\label{Py_ptwise}
Let $F_n, F \colon [0, +\infty) \to [0, +\infty)$, $n = 1, 2, \ldots$, be continuous metric preserving functions. Assume that, for any mm-space $X$, the metric transformed space $F_n(X)$ concentrates to $F(X)$ as $n \to \infty$. Then, $F_n$ converges pointwise to $F$ as $n \to \infty$.
\end{prop}

\begin{proof}
Take any real number $s > 0$ and fix it. An mm-space $X$ is defined as
\begin{equation*}
X := (\{0, s\}, |\cdot|, \delta_{0, s}),
\end{equation*}
where $\delta_{0, s} := 1/2 \delta_0 + 1/2 \delta_s$ and $\delta_x$ is the Dirac probability measure at $x$. By the assumption, $F_n(X)$ concentrates to $F(X)$ as $n \to \infty$. Note that $F_n(X)$ is mm-isomorphic to the mm-space
\begin{equation*}
(\{0, F_n(s)\}, |\cdot|, \delta_{0, F_n(s)}).
\end{equation*}
Suppose that $F_n(s)$ diverges to infinity as $n \to \infty$. It is easy to see that $F_n(X)$ converges weakly to the pyramid
\begin{equation*}
\Py := \left\{(\{0, t\}, |\cdot|, \delta_{0, t}) \midd t \geq 0\right\}.
\end{equation*}
Since $\Py$ is not $\square$-precompact (see \cite{MMG}*{Lemma 4.28}), it follows from Lemma \ref{chara_ass_py} that $\Py \neq \Py{F(X)}$, which is a contradiction. Thus we can assume that $F_n(s)$ converges to a nonnegative number $t$. Then, $F_n(X)$ $\square$-converges to $(\{0, t\}, |\cdot|, \delta_{0, t})$, so that we have $t = F(s)$. This completes the proof.
\end{proof}

The following lemma was obtained in \cite{prod}.

\begin{lem}[\cite{prod}*{Claim 5.1}]\label{prod:claim5.1}
Let $F_n, F \colon [0, +\infty) \to [0, +\infty)$, $n = 1, 2, \ldots$, be continuous metric preserving functions. Assume that $F_n$ converges pointwise to $F$ as $n \to \infty$. If \eqref{liminf:eq} does not hold, then there exist a sequence $\{X_n\}_{n \in \N}$ of mm-spaces and two mm-spaces $X$, $Y$ such that 
\begin{enumerate}
\item $X_n$ concentrates to $X$ as $n \to \infty$,
\item $F_n(X_n)$ concentrates to $Y$ as $n \to \infty$,
\item $F(X)$ and $Y$ are not mm-isomorphic to each other.
\end{enumerate}
Moreover, if $F$ is increasing, then there exists an mm-space $X'$ such that $Y = F(X')$.
\end{lem}

\begin{prop}\label{limsup_ness}
Let $F_n, F \colon [0, +\infty) \to [0, +\infty)$, $n = 1, 2, \ldots$, be continuous metric preserving functions. Assume that $F_n$ converges pointwise to $F$ as $n \to \infty$ and that $F$ is nondecreasing. If
\begin{equation*}
\limsup_{n \to \infty} \sup{F_n} > \sup{F},
\end{equation*}
then there exists a sequence $\{X_n\}_{n \in \N}$ of mm-spaces converging weakly to a pyramid $\Py$ such that $F_n(X_n)$ does not converge weakly to $F(\Py)$.
\end{prop}

\begin{proof}
Assume that $\limsup_{n \to \infty} \sup{F_n} > \sup{F}$. We define a pyramid
\begin{equation*}
\Py := \left\{(\{0, s\}, |\cdot|, \delta_{0, s}) \midd s \geq 0\right\} \subset \X,
\end{equation*}
where the notation is same as in the proof of Proposition \ref{Py_ptwise}. We set
\begin{equation*}
\alpha := \limsup_{n\to\infty} \sup{F_n} \quad \text{ and } \quad \eta := \left\{\begin{array}{ll} \frac{\alpha -\sup{F}}{2} & \text{ if } \alpha < +\infty, \\ 1 & \text{ if } \alpha = +\infty. \end{array}\right.
\end{equation*}
Note that $\eta > 0$. There exists a subsequence $\{n_i\} \subset \{n\}$ such that
\begin{equation*}
\sup{F_{n_i}} > \sup{F} + \eta.
\end{equation*}
For each $i$, there exists $s_i \in [0, +\infty)$ such that
\begin{equation*}
F_{n_i}(s_i) > \sup{F} + \eta.
\end{equation*}
We see that $s_i \to \infty$ as $i \to \infty$. In fact, if
\begin{equation*}
s := \liminf_{i \to \infty} s_i < +\infty,
\end{equation*}
then, by Lemma \ref{ptcpt}, we have
\begin{equation*}
F(s) \geq \sup{F} + \eta,
\end{equation*}
which is a contradiction. We define mm-spaces $X_i$, $i = 1, 2, \ldots$, by 
\begin{equation*}
X_i := (\{0, s_i\}, |\cdot|, \delta_{0, s_i}),
\end{equation*}
and then it follows from $s_i \to \infty$ that $X_i$ converges weakly to $\Py$ as $i \to \infty$. We prove that $F_{n_i}(X_i)$ does not converge weakly to $F(\Py)$. Since $\sup{F} + \eta < F_{n_i}(s_i)$ for any $i$, we have
\begin{equation*}
(\{0, \sup{F} + \eta\}, |\cdot|, \delta_{0, \sup{F} + \eta}) \in \Py {F_{n_i}(X_i)},
\end{equation*}
which implies
\begin{equation*}
\liminf_{n \to \infty} \square((\{0, \sup{F} + \eta\}, |\cdot|, \delta_{0, \sup{F} + \eta}), \Py {F_{n_i}(X_i)}) = 0.
\end{equation*}
On the other hand, since $\diam{Y} \leq \sup{F}$ for any $Y \in F(\Py)$, we have
\begin{equation*}
(\{0, \sup{F} + \eta\}, |\cdot|, \delta_{0, \sup{F} + \eta}) \not\in F(\Py).
\end{equation*}
Thus we obtain the conclusion.
\end{proof}

\begin{proof}[Proof of {\rm`(1) $\Rightarrow$ (2)'} of Theorem \ref{Py}]
This follows from Proposition \ref{Py_ptwise}, Lemma \ref{prod:claim5.1}, and Proposition \ref{limsup_ness} directly.
\end{proof}

We prepare some lemmas for the proof of  {\rm`(2) $\Rightarrow$ (1)'} of Theorem \ref{Py}.

\begin{lem}\label{lip_up}
Let $F_n, F \colon [0,+\infty) \to [0,+\infty)$, $n= 1, 2, \ldots$, be continuous metric preserving functions satisfying \eqref{ptwise} and \eqref{liminf:eq}. If Borel measurable maps $f_n \colon X_n \to X$ between mm-spaces $X_n$ and $X$ are $1$-Lipschitz up to $\varepsilon_n$ for some $\varepsilon_n \to 0$, and if $\prok^{d_X}({f_n}_* m_{X_n}, m_X) \leq \varepsilon_n$ holds, then there exists a sequence $\delta_n \to 0$ such that $f_n$ is $1$-Lipschitz up to $\delta_n$ with respect to $F_n \circ d_{X_n}$ and $F \circ d_X$, and
\begin{equation*}
\prok^{F \circ d_X}({f_n}_* m_{X_n}, m_X) \leq \delta_n.
\end{equation*}
\end{lem}

\begin{proof}
We take any real number $\varepsilon > 0$. By the inner regularity of $m_X$, there exists a compact set $K \subset X$ such that $m_X(K) \geq 1 - \varepsilon$. We put $D_\varepsilon := \diam{K} + 2\varepsilon$. By the assumptions of this lemma, Lemma \ref{ptcpt}, and Lemma \ref{liminf:lem}, for any sufficiently large $n$,
\begin{itemize}
\item $f_n \colon X_n \to X$ is $1$-Lipschitz up to $\varepsilon$ and $\prok({f_n}_* m_{X_n}, m_X) \leq \varepsilon$ holds,
\item $|F_n(s) - F(s)| < \varepsilon$ holds for any $s \in [0, D_\varepsilon]$,
\item $F_n(s) \leq F_n(s') + \varepsilon$ holds for any $s \in [0, D_\varepsilon]$ and for any $s' \geq s$.
\end{itemize}
Let $X'_n \subset X_n$ be a nonexceptional domain of $f_n \colon X_n \to X$ and let
\begin{equation*}
\widetilde{X}_n := X'_n \cap f_n^{-1}(U_\varepsilon(K)).
\end{equation*}
We see that
\begin{equation*}
m_{X_n}(\tilde{X_n}) \geq m_{X_n}(X'_n) + {f_n}_* m_{X_n}(U_\varepsilon(K)) - 1 \geq 1 - 3\varepsilon.
\end{equation*}
For any $x, x' \in \widetilde{X}_n$, we have
\begin{align*}
&F(d_X(f_n(x), f_n(x'))) \leq F_n(d_X(f_n(x), f_n(x'))) + \varepsilon  \\
&\leq F_n(d_{X_n}(x, x') + \varepsilon) + 2\varepsilon \leq F_n(d_{X_n}(x, x')) + F_n(\varepsilon) + 2\varepsilon \\
&\leq F_n(d_{X_n}(x, x')) + F(\varepsilon) + 3\varepsilon,
\end{align*}
where the first and second inequalities follow from 
\begin{equation*}
\diam{f_n(\widetilde{X}_n)} \leq D_\varepsilon.
\end{equation*}
Thus, for any sufficiently large $n$, the map $f_n \colon F_n(X_n) \to F(X)$ is $1$-Lipschitz up to $F(\varepsilon) + 3\varepsilon$. 

We next prove $\prok^{F \circ d_X}({f_n}_* m_{X_n}, m_X) \leq \max{\{\varepsilon, F(\varepsilon)\}}$. For any subset $A \subset X$ and any real number $\eta > 0$, we have
\begin{equation*}
U_\varepsilon^{d_X}(A) \subset U_{F(\varepsilon) + \eta}^{F \circ d_X}(A).
\end{equation*}
In fact, taking any point $x \in U_\varepsilon^{d_X}(A)$, it holds that $d_X(x, A) < \varepsilon$, which implies that $F(d_X(x, A)) \leq F(\varepsilon)$. Combining $\prok^{d_X}({f_n}_* m_{X_n}, m_X) \leq \varepsilon$ and this leads to
\begin{equation*}
m_X(A) \leq {f_n}_* m_{X_n}(U_\varepsilon^{d_X}(A)) + \varepsilon \leq {f_n}_* m_{X_n}(U_{F(\varepsilon) + \eta}^{F \circ d_X}(A)) + \varepsilon,
\end{equation*}
which implies that $\prok^{F \circ d_X}({f_n}_* m_{X_n}, m_X) \leq \max{\{\varepsilon, F(\varepsilon) + \eta\}}$. As $\eta \to 0$, we obtain 
\begin{equation*}
\prok^{F \circ d_X}({f_n}_* m_{X_n}, m_X) \leq \max{\{\varepsilon, F(\varepsilon)\}}.
\end{equation*}
The proof is completed.
\end{proof}

\begin{lem}\label{Py:lem} 
Let $X$ and $Y$ be two mm-spaces and let $\varepsilon > 0$.
\begin{enumerate}
\item If a Borel measurable map $f \colon X \to Y$ is $1$-Lipschitz up to $\varepsilon$ and if $\prok(f_* m_X, m_Y) \leq \varepsilon$ holds, then $\square(Y, \Py X) \leq 4\varepsilon$.
\item If $\square(Y, \Py X) < \varepsilon$, then there exists a Borel measurable map $f \colon X \to Y$ that is $1$-Lipschitz up to $3\varepsilon$ and $\prok(f_* m_X, m_Y) \leq 3\varepsilon$ holds.
\end{enumerate}
\end{lem}

\begin{proof}
We first prove (1). Let $f \colon X \to Y$ be a map satisfying the assumption. By Lemma \ref{fin_dim_app}, there exist $1$-Lipschitz maps $\Phi_N \colon Y \to (\R^N, \|\cdot\|_\infty)$, $N=1, 2, \ldots$, such that
\begin{equation*}
\lim_{N\to\infty} \square(Y, (\R^N, \|\cdot\|_{\infty}, {\Phi_N}_*m_Y)) = 0.
\end{equation*}
Since the composition $\Phi_N \circ f \colon X \to (\R^N, \|\cdot\|_{\infty})$ is $1$-Lipschitz up to $\varepsilon$, by Lemma \ref{MWext}, there exists a $1$-Lipschitz map $\Psi_N \colon X \to (\R^N, \|\cdot\|_\infty)$ such that
\begin{equation*}
\kf^{m_X}(\Phi_N \circ f, \Psi_N) \leq \varepsilon.
\end{equation*}
Note that $(\R^N, \|\cdot\|_\infty, {\Psi_N}_* m_X) \in \Py X$. Then we have 
\begin{align*}
& \prok({\Phi_N}_* m_Y, {\Psi_N}_* m_X) \\
& \leq \prok((\Phi_N \circ f)_* m_X, {\Psi_N}_* m_X) + \prok((\Phi_N \circ f)_* m_X, {\Phi_N}_* m_Y)  \\
& \leq \kf^{m_X}(\Phi_N \circ f, \Psi_N) + \prok(f_* m_X, m_Y)  \leq 2\varepsilon,
\end{align*}
where the second inequality follows from Lemma \ref{prok_kf} and \eqref{lip_prok}. Thus we have
\begin{align*}
&\square(Y, \Py X) \leq \square(Y, (\R^N, \|\cdot\|_\infty, {\Psi_N}_* m_X)) \\
&\leq \square(Y, (\R^N, \|\cdot\|_\infty, {\Phi_N}_* m_{Y})) + 2\prok({\Phi_N}_* m_Y, {\Psi_N}_* m_X) \\
&\leq \square(Y, (\R^N, \|\cdot\|_\infty, {\Phi_N}_* m_{Y})) + 4\varepsilon.
\end{align*}
As $N \to \infty$, we obtain $\square(Y, \Py X) \leq 4\varepsilon$.

We next prove (2). There exists an mm-space $X'$ dominated by $X$ such that $\square(Y, X') < \varepsilon$. There exist a $1$-Lipschitz map $f \colon X \to X'$ with $f_* m_X = m_{X'}$ and a $3\varepsilon$-mm-isomorphism $g \colon X' \to Y$ by Lemma \ref{lem4.22}. It is easy to see that the composition $g \circ f \colon X \to Y$ is $1$-Lipschitz up to $3\varepsilon$ and that $\prok((g \circ f)_* m_X, m_Y) \leq 3\varepsilon$. This completes the proof.
\end{proof}

We see that there exists an mm-space $X' \in \Py X$ minimizing $\square(Y, \Py X)$ by Lemma \ref{chara_ass_py}.

\begin{cor}\label{Py:lem_cor}
Let $\Py$ be a pyramid and let $Y$ an mm-space. Assume that, for any $\varepsilon > 0$, there exist an mm-space $X_\varepsilon \in \Py$ and a Borel measurable map $f_\varepsilon \colon X_\varepsilon \to Y$ such that $f_\varepsilon$ is $1$-Lipschitz up to $\varepsilon$ and $\prok({f_\varepsilon}_* m_{X_\varepsilon}, m_Y) \leq \varepsilon$ holds. Then $Y \in \Py$.
\end{cor}

\begin{proof}
Take any $\varepsilon > 0$. There exist an mm-space $X_\varepsilon \in \Py$ and a map $f_\varepsilon \colon X_\varepsilon \to Y$ in the assumption. By Lemma \ref{Py:lem}, we have $\square(Y, \Py X_\varepsilon) \leq 4\varepsilon$. Since $\Py X_\varepsilon \subset \Py$, we have
\begin{equation*}
\square(Y, \Py) \leq \square(Y, \Py X_\varepsilon) \leq 4\varepsilon.
\end{equation*}
As $\varepsilon \to 0$, we obtain $Y \in \Py$.
\end{proof}

We prove {\rm`(2) $\Rightarrow$ (1)'} of Theorem \ref{Py}. Let $F_n$, $F$ be continuous metric preserving functions satisfying \eqref{ptwise}, \eqref{liminf:eq}, and \eqref{limsup:eq}. Take any sequence $\{X_n\}$ of mm-spaces and any pyramid $\Py \in \Pi$ such that $X_n$ converges weakly to $\Py$. Our goal is to prove that $F_n(X_n)$ converges weakly to $F(\Py)$ as $n \to \infty$.

\begin{prop}\label{liminf:prop}
For any $Y \in F(\Py)$, we have 
\begin{equation*}
\lim_{n \to \infty}\square(Y, \Py{F_n(X_n)}) = 0.
\end{equation*}
\end{prop}

\begin{proof}[Proof of Proposition \ref{liminf:prop}]
Take any mm-space $Y \in F(\Py)$ and any real number $\varepsilon > 0$. There exist two mm-spaces $X' \in \Py$ and $Y' \in \X$ such that $F(X')$ dominates $Y'$ and $\square(Y', Y) \leq \varepsilon$. Since $X_n$ converges weakly to $\Py$, we have $\square(X', \Py X_n) \to 0$ as $n \to \infty$. Then, by Lemma \ref{Py:lem}, there exist Borel measurable maps $f_n \colon X_n \to X'$, $n = 1, 2, \ldots$, such that $f_n$ is $1$-Lipschitz up to $\varepsilon_n$ and $\prok({f_n}_* m_{X_n}, m_{X'}) \leq \varepsilon_n$ holds for some $\varepsilon_n \to 0$. By Lemma \ref{lip_up}, we see that the map $f_n \colon F_n(X_n) \to F(X')$ is $1$-Lipschitz up to $\delta_n$ and $\prok^{F \circ d_{X'}}({f_n}_* m_{X_n}, m_{X'}) \leq \delta_n$ holds for some $\delta_n \to 0$. Moreover, since $F(X')$ dominates $Y'$, there exists a $1$-Lipschitz map $g \colon F(X') \to Y'$ with $g_* m_{X'} = m_{Y'}$. Since the composition $g \circ f_n \colon F_n(X_n) \to Y'$ is also $1$-Lipschitz up to $\delta_n$ and fulfills $\prok({(g \circ f_n)}_* m_{X_n}, m_{Y'}) \leq \delta_n$, we have
\begin{equation*}
\square(Y', \Py{F_n(X_n)}) \leq 4\delta_n
\end{equation*}
by Lemma \ref{Py:lem}. Thus we have
\begin{equation*}
\square(Y, \Py{F_n(X_n)}) \leq \square(Y, Y') + \square(Y', \Py{F_n(X_n)}) \leq \varepsilon + 4\delta_n.
\end{equation*}
As $n \to \infty$ and $\varepsilon \to 0$, we obtain the conclusion.
\end{proof}

\begin{prop}\label{limsup:prop}
If an mm-space $Y$ satisfies
\begin{equation*}
\liminf_{n \to \infty} \square(Y, \Py F_n(X_n)) = 0,
\end{equation*}
then $Y \in F(\Py)$.
\end{prop}

The following proof is inspired by the ideas invented first by Ryunosuke Ozawa and Takumi Yokota. The author was  privately informed of their ideas.

\begin{proof}
Assume that $\liminf_{n \to \infty} \square(Y, \Py F_n(X_n)) = 0$. We can assume that $Y \neq \{*\}$. Choosing a subsequence of $n$, we can assume that $\square(Y, \Py F_n(X_n)) \to 0$ as $n \to \infty$. Then, by Lemma \ref{Py:lem}, there exist Borel measurable maps $f_n \colon F_n(X_n) \to Y$ and a sequence $\varepsilon_n \to 0$ such that $f_n$ is $1$-Lipschitz up to $\varepsilon_n$ and $\prok({f_n}_* m_{X_n}, m_Y) \leq \varepsilon_n$ holds for every $n$. Let $\widetilde{X}_n \subset X_n$ be a nonexceptional domain of $f_n \colon F_n(X_n) \to Y$.

We take any sufficiently small real number $\varepsilon > 0$. We find finite many open sets $B_1, \ldots, B_N$ in $Y$ such that $\diam{B_i} < \varepsilon$ and $m_Y(B_i) > 0$ for every $i \in \{1, \ldots, N\}$ and that
\begin{equation*}
\delta' := \min_{1 \leq i < j \leq N} d_Y(B_i, B_j)  > 0, \quad \sum_{i=1}^N m_Y(B_i)  > 1 - \varepsilon.
\end{equation*}
Let $B_0 := Y \setminus \bigsqcup_{i=1}^N B_i$. For any $i = 0, \ldots, N$, we take any point $y_i \in B_i$ and fix it. If $B_0 = \emptyset$, we consider only for $i = 1, \ldots, N$. An mm-space $\dot{Y}$ is defined as
\begin{equation*}
\dot{Y} := (\{y_i\}_{i=0}^N, d_Y, m_{\dot{Y}}),
\end{equation*}
where $m_{\dot{Y}}(\{y_i\}) := m_Y(B_i)$. Note that the natural embedding $\iota \colon \dot{Y} \ni y_i \mapsto y_i \in Y$ is an $\varepsilon$-mm-isomorphism. Our goal is, by Corollary \ref{Py:lem_cor}, to prove that there exist an mm-space $W_\varepsilon \in F(\Py)$ and a map $h_\varepsilon \colon W_\varepsilon \to Y$ such that $h_\varepsilon$ is $1$-Lipschitz up to $4\varepsilon$ and $\prok({h_\varepsilon}_* m_{W_\varepsilon}, m_Y) \leq 2\varepsilon$.

For any $n \in \N$ and for any $i = 1, \ldots, N$, we define a subset 
\begin{equation*}
A_{n, i}  := f_n^{-1}(B_i) \cap \widetilde{X}_n \subset X_n
\end{equation*}
and define a real number
\begin{equation*}
R := \inf\left\{r > 0 \midd \max_{i,j=1, \ldots, N}{d_Y(y_i, y_j)} -3\varepsilon \leq F(r) \right\}.
\end{equation*}
The existence of $R$ is discussed as follows. For any $i, j = 1, \ldots, N$, taking any $x_{n, i} \in A_{n, i}$ and any $x_{n, j} \in A_{n, j}$, we have
\begin{equation*}
\begin{split}
d_Y(y_i, y_j) & \leq d_Y(f_n(x_{n,i}), f_n(x_{n,j})) + 2\varepsilon \\
& \leq F_n(d_{X_n}(x_{n,i}, x_{n,j})) + \varepsilon_n + 2\varepsilon \leq \sup{F_n} + \varepsilon_n + 2\varepsilon.
\end{split}
\end{equation*}
By \eqref{limsup:eq}, we obtain $d_Y(y_i, y_j) - 2\varepsilon \leq \sup{F}$, which means that $R$ exists. 

We define a map $\Phi_n \colon X_n \to (\R^N, \|\cdot\|_{\infty})$ by
\begin{equation*}
\Phi_n(x) : = (\min\{d_{X_n}(x, A_{n, i}), R\})_{i=1}^N, \quad x \in X_n.
\end{equation*}
Since the measure $\nu_n := {\Phi_n}_* m_{X_n}$ is supported in the compact set $\left\{z \in \R^N \midd  \|z\|_{\infty} \leq R \right\}$, a sequence $\{\nu_n\}_{n \in \N}$ is tight. Thus the sequence $\{\nu_n\}$ has a subsequence converging weakly to a Borel probability measure $\nu$ on $\R^N$. An mm-space $Z$ is defined as
\begin{equation*}
Z : = (\supp{\nu}, \|\cdot\|_{\infty}, \nu).
\end{equation*}
It follows from $(\R^N, \|\cdot\|_{\infty}, \nu_n) \prec X_n$ that $\square(Z, \Py{X_n}) \to 0$ as $n \to \infty$. Since $X_n$ converges weakly to $\Py$, we have $Z \in \Py$.

Let
\begin{equation*}
\delta := \frac{1}{2} \inf\left\{r>0 \midd \frac{\delta'}{2} \leq F(r)\right\} > 0
\end{equation*}
and let 
\begin{equation*}
A_i := \left\{(z_i)_{i=1}^N \in \R^N \midd z_i = 0 \text{ and } |z_j| \geq \delta \text{ for any } j \neq i\right\}
\end{equation*}
for any $i=1, \ldots, N$. Note that $2\delta \leq R$. In fact, if $2\delta > R$, then we have
\begin{equation*}
\frac{\delta'}{2} > F(R) \geq \max_{i,j =1, \ldots, N}{d_Y(y_i, y_j)} - 3\varepsilon \geq \delta' - 3\varepsilon,
\end{equation*}
which implies that $\max_{i,j}{d_Y(y_i, y_j)}  \leq 6\varepsilon$. This contradicts $Y \neq \{*\}$.

\begin{claim}\label{claim_inc}
For any $i =1, \ldots, N$ and for any sufficiently large $n$, we have
\begin{equation*}
\Phi_n(A_{n, i}) \subset A_i.
\end{equation*}
\end{claim}

\begin{proof}
Take any $i \in \{1, \ldots, N\}$ and fix it. Since
\begin{equation*}
(\Phi_n(x))_i = \min{\{d_{X_n}(x, A_{n, i}), R\}} = 0
\end{equation*}
for any $x \in A_{n, i}$, it is sufficient to prove that  
\begin{equation}\label{claim_inc:eq}
\liminf_{n \to \infty} d_{X_n}(A_{n,i}, A_{n, j}) \geq 2\delta
\end{equation}
for every $j \neq i$. In fact, taking $2\delta \leq R$ into account, we have
\begin{equation*}
\liminf_{n \to \infty} \inf_{x \in A_{n, i}} (\Phi_n(x))_j \geq 2\delta,
\end{equation*}
which implies that $\inf_{x \in A_{n, i}} (\Phi_n(x))_j \geq \delta$ for any sufficiently large $n$. We prove \eqref{claim_inc:eq}. We put $\alpha := \liminf_{n \to \infty} d_{X_n}(A_{n,i}, A_{n, j})$. There exists a subsequence $\{n_k\} \subset \{n\}$ such that $d_{X_{n_k}}(A_{n_k,i}, A_{n_k, j}) \to \alpha$ as $k \to \infty$. We find two sequences $\{x_m\}_m \subset A_{n_k, i}$ and $\{x'_m\}_m \subset A_{n_k, j}$ such that $d_{X_{n_k}}(x_m, x'_m) \to d_{X_{n_k}}(A_{n_k,i}, A_{n_k, j})$ as $m \to \infty$. Then we have
\begin{align*}
\delta' &\leq d_Y(B_i, B_j) \leq d_{Y}(f_n(x_m), f_n(x'_m)) \\
&\leq F_{n_k}(d_{X_{n_k}}(x_m, x'_m)) + \varepsilon_{n_k} \\
& \to  F_{n_k}(d_{X_{n_k}}(A_{n_k, i}, A_{n_k, j})) + \varepsilon_{n_k} \text{ as }  m \to \infty
\end{align*}
for each $k$. Moreover, by Lemma \ref{ptcpt}, we have
\begin{equation*}
\lim_{k \to \infty} F_{n_k}(d_{X_{n_k}}(A_{n_k, i}, A_{n_k, j})) = F(\alpha).
\end{equation*}
Combining these implies $\delta' \leq F(\alpha)$, so that $2\delta \leq \alpha$. Therefore we obtain \eqref{claim_inc:eq}. This completes the proof.
\end{proof}

We define a subset $\widetilde{Z} \subset Z$ by
\begin{equation*}
\widetilde{Z} := \supp{\nu} \cap \bigsqcup_{i=1}^N A_i
\end{equation*}
and define a map $g \colon F(Z) \to \dot{Y}$ by
\begin{equation*}
g(z) := \left\{ \begin{array}{ll} y_i & \text{ if } z \in A_i, \\ y_0 & \text{ if } z \in Z \setminus \widetilde{Z}. \end{array}\right.
\end{equation*}
If $B_0 = \emptyset$, then we set $g(z)$, for $z \in Z \setminus \widetilde{Z}$, an arbitrary point of $\dot{Y}$.

\begin{claim}\label{claim_1-Lip}
The map $g \colon F(Z) \to \dot{Y}$ is $1$-Lipschitz up to $3\varepsilon$ and fulfills
\begin{equation*}
\prok(g_* \nu, m_{\dot{Y}}) \leq \varepsilon.
\end{equation*}
\end{claim}

If we prove this claim, then the composition $\iota \circ g \colon F(Z) \to Y$ is $1$-Lipschitz up to $4\varepsilon$ and $\prok((\iota \circ g)_* \nu, m_Y) \leq 2\varepsilon$ holds. Combining this and Corollary \ref{Py:lem_cor} implies $Y \in F(\Py)$. 

\begin{proof}[Proof of Claim \ref{claim_1-Lip}]
We first prove $\prok(g_* \nu, m_{\dot{Y}}) \leq \varepsilon$. For any $i = 1, \ldots, N$, we have
\begin{align*}
m_{\dot{Y}}(\{y_i\}) &= m_Y(B_i) \leq \liminf_{n \to \infty} {f_n}_* m_{X_n}(B_i) = \liminf_{n \to \infty} m_{X_n}(A_{n, i})\\
&\leq \liminf_{n \to \infty} {\Phi_n}_* m_{X_n}(A_i) \leq \nu(A_i) = g_* \nu(\{y_i\}),
\end{align*}
where the second inequality follows from Claim \ref{claim_inc}. Moreover, since $m_{\dot{Y}}(\{y_0\}) \leq \varepsilon$, we have
\begin{equation*}
m_{\dot{Y}}(B) \leq g_* \nu(B) + \varepsilon
\end{equation*}
for any subset $B \subset \dot{Y}$, which implies that $\prok(g_* \nu, m_{\dot{Y}}) \leq \varepsilon$.

We next prove that the map $g \colon F(Z) \to \dot{Y}$ is $1$-Lipschitz up to $3\varepsilon$. Since
\begin{equation*}
\nu(Z \setminus \widetilde{Z}) = g_* \nu(\{y_0\}) \leq m_{\dot{Y}}(\{y_0\}) \leq \varepsilon,
\end{equation*}
we see that $\nu(\widetilde{Z}) \geq 1- \varepsilon$. It is sufficient to prove that 
\begin{equation}\label{claim_1-Lip:eq}
d_{\dot{Y}}(g(z), g(z')) \leq F(\|z - z'\|_\infty) + 3\varepsilon
\end{equation}
holds for any $z, z' \in \widetilde{Z}$. We take any $i, j \in \{1, \ldots, N\}$ and fix them. We take any points $z \in A_i$ and $z' \in A_j$. There exist $x_n \in A_{n, i}$ and $x'_n \in A_{n, j}$, $n =1, 2, \ldots$, such that 
\begin{equation*}
\|\Phi_n(x_n) - z \|_\infty, \|\Phi_n(x'_n) - z' \|_\infty \to 0 \text{ as } n \to \infty.
\end{equation*}
Since we have
\begin{equation*}
d_Y(y_i, y_j) \leq d_Y(B_i, B_j) + 2\varepsilon \leq d_Y(f_n(x_n), f_n(\tilde{x}')) + 2\varepsilon \leq F_n(d_{X_n}(x_n, \tilde{x}')) + \varepsilon_n + 2\varepsilon
\end{equation*}
for any $\tilde{x}' \in A_{n, j}$, we obtain
\begin{equation*}
d_Y(y_i, y_j) \leq \limsup_{n\to\infty} F_n(d_{X_n}(x_n, A_{n, j})) + 2\varepsilon.
\end{equation*}
On the other hand, we see that
\begin{equation*}
d_Y(y_i, y_j) \leq \max_{i, j=1, \ldots, N} d_Y(y_i, y_j) \leq F(R) + 3\varepsilon = \lim_{n \to \infty} F_n(R) + 3\varepsilon.
\end{equation*}
Combining these implies that
\begin{align*}
d_Y(y_i, y_j) &\leq \limsup_{n\to\infty} F_n(\min\{d_{X_n}(x_n, A_{n, j}), R\}) + 3\varepsilon \\
& =  \limsup_{n\to\infty} F_n((\Phi_n(x_n))_j) + 3\varepsilon \\
& = \limsup_{n\to\infty} F_n((\Phi_n(x_n) - \Phi_n(x'_n))_j) + 3\varepsilon.
\end{align*}
By \eqref{liminf:eq} and Lemma \ref{liminf:lem}, we have 
\begin{equation*}
d_Y(y_i, y_j) \leq \limsup_{n\to\infty} F_n(\|\Phi_n(x_n) - \Phi_n(x'_n)\|_\infty) + 3\varepsilon.
\end{equation*}
Furthermore, since  $\|\Phi_n(x_n) - \Phi_n(x'_n)\|_\infty \to \|z - z'\|_\infty$ as $n \to \infty$, we have
\begin{equation*}
d_{\dot{Y}}(g(z), g(z')) = d_Y(y_i, y_j)  \leq F(\|z - z'\|_\infty) + 3\varepsilon
\end{equation*}
by Lemma \ref{ptcpt}. Therefore we obtain \eqref{claim_1-Lip:eq}. This completes the proof.
\end{proof}

Applying Corollary \ref{Py:lem_cor} to the map $\iota \circ g \colon F(Z) \to Y$, we obtain $Y \in F(\Py)$. This completes the proof of Proposition \ref{limsup:prop}.
\end{proof}

\begin{proof}[Proof of {\rm`(2) $\Rightarrow$ (1)'} of Theorem \ref{Py}]
This follows from Proposition \ref{liminf:prop} and Proposition \ref{limsup:prop}.
\end{proof}

\begin{proof}[Proof of Corollary \ref{Py:cor} {\rm (A)}]
`(A1) $\Rightarrow$ (A2)' follows from Proposition \ref{Py_ptwise} and Proposition \ref{limsup_ness}. 

We prove `(A2) $\Rightarrow$ (A1)'. Assume that $\Py_n$ converges weakly to $\Py$.

We take any mm-space $Y \in F(\Py)$ and any real number $\varepsilon > 0$. There exist two mm-spaces $X \in \Py$ and $Y' \in \X$ such that $F(X)$ dominates $Y'$ and $\square(Y, Y') \leq \varepsilon$. Since $\Py_n$ converges weakly to $\Py$, there exist mm-spaces $X_n \in \Py_n$, $n = 1, 2, \ldots$, such that $\square(X_n, X) \to 0$ as $n \to \infty$. Since, in particular, $X_n$ converges weakly to $\Py X$ and $Y' \in \Py F(X)$, we have
\begin{equation*}
\lim_{n \to \infty} \square(Y', \Py F_n(X_n)) = 0
\end{equation*}
by Proposition \ref{liminf:prop}. Thus we have
\begin{equation*}
\square(Y, F_n(\Py_n)) \leq \square(Y', F_n(\Py_n)) + \varepsilon \leq \square(Y', \Py F_n(X_n)) + \varepsilon,
\end{equation*}
which implies that $\lim_{n \to \infty} \square(Y, F_n(\Py_n)) = 0$.

On the other hand, we assume that an mm-space $Y$ satisfies
\begin{equation*}
\liminf_{n \to \infty} \square(Y, F_n(\Py_n)) = 0.
\end{equation*}
Then, there exist mm-spaces $X_n \in \Py_n$, $n = 1, 2, \ldots$, such that
\begin{equation*}
\liminf_{n \to \infty} \square(Y, \Py F_n(X_n)) = 0.
\end{equation*}
By the compactness of $\Pi$, $\{X_n\}_{n \in \N}$ has a weak convergent subsequence. Let $\PyQ$ be a limit pyramid of a weak convergent subsequence of $\{X_n\}_{n \in \N}$. By Proposition \ref{limsup:prop}, we have $Y \in F(\PyQ)$. Moreover, since $X_n \in \Py_n$ for each $n$, we have $\PyQ \subset \Py$, which implies that $F(\PyQ) \subset F(\Py)$. Thus we obtain $Y \in F(\Py)$.

Combining these means that $F_n(\Py_n)$ converges weakly to $F(\Py)$. This completes the proof of this corollary.
\end{proof}

\subsection{Proof of Theorem \ref{Py_con}}

We start with proving the following lemma. 

\begin{lem}\label{lip_up_con}
Let $F_n, F \colon [0,+\infty) \to [0,+\infty)$, $n= 1, 2, \ldots$, be continuous metric preserving functions satisfying \eqref{ptwise}, \eqref{liminf:eq} and \eqref{incr}. If
\begin{equation*}
\lim_{n \to \infty}\square(F(Y), \Py{F_n(X_n)}) = 0
\end{equation*} 
for mm-spaces $Y$ and $X_n$, $n = 1, 2, \ldots$, then we have 
\begin{equation*}
\lim_{n \to \infty}\square(Y, \Py{X_n}) = 0.
\end{equation*}
\end{lem}

\begin{proof}
Assume that $\square(F(Y), \Py{F_n(X_n)}) \to 0$ as $n \to \infty$. By Lemma \ref{Py:lem}, there exist Borel measurable maps $f_n \colon F_n(X_n) \to F(Y)$ and a sequence $\varepsilon_n \to 0$ such that $f_n$ is $1$-Lipschitz up to $\varepsilon_n$ and fulfills $\prok^{F\circ d_Y}({f_n}_* m_{X_n}, m_Y) \leq \varepsilon_n$ for every $n$.

We take any real number $\varepsilon > 0$. We prove the following: there exist an mm-space $Z$  and a map $h \colon Z \to Y$ (both depending on $\varepsilon$) such that 
\begin{itemize}
\item $\square(Z, \Py{X_n}) \to 0$ as $n \to \infty$, 
\item $h$ is $1$-Lipschitz up to $3\varepsilon$ and $\prok(h_* m_Z, m_Y) \leq 2\varepsilon$ holds.
\end{itemize}
If we prove this, then, by Lemma \ref{Py:lem}, we have
\begin{equation*}
\limsup_{n \to \infty} \square(Y, \Py{X_n}) \leq 12\varepsilon.
\end{equation*}

Note that the outline of the following proof is similar to the proof of Proposition \ref{limsup:prop}.

We find finite many open sets $B_1, \ldots, B_N \subset Y$ such that $\diam{B_i} < \varepsilon$ and $m_Y(B_i) > 0$ for every $i \in \{1, \ldots, N\}$ and that
\begin{equation*}
\delta := \min_{1 \leq i < j \leq N} d_Y(B_i, B_j)  > 0, \quad \sum_{i=1}^N m_Y(B_i)  > 1 - \varepsilon.
\end{equation*}
Let $B_0 := Y \setminus \bigsqcup_{i=1}^N B_i$. For any $i = 0, \ldots, N$, we take any point $y_i \in B_i$ and fix it. If $B_0 = \emptyset$, we consider only for $i = 1, \ldots, N$. An mm-space $\dot{Y}$ is defined as
\begin{equation*}
\dot{Y} := (\{y_i\}_{i=0}^N, d_Y, m_{\dot{Y}}),
\end{equation*}
where $m_{\dot{Y}}(\{y_i\}) := m_Y(B_i)$. Note that the natural embedding $\iota \colon \dot{Y} \ni y_i \mapsto y_i \in Y$ is an $\varepsilon$-mm-isomorphism.

For any $n \in \N$ and for any $i = 1, \ldots, N$, we define a subset 
\begin{equation*}
A_{n, i}  := f_n^{-1}(B_i) \cap \widetilde{X}_n \subset X_n
\end{equation*}
and define a real number
\begin{equation*}
R := \max_{i,j=1, \ldots, N}{d_Y(y_i, y_j)} .
\end{equation*}
We see that $\delta \leq R$. We define a map $\Phi_n \colon X_n \to (\R^N, \|\cdot\|_{\infty})$ by
\begin{equation*}
\Phi_n(x) : = (\min\{d_{X_n}(x, A_{n, i}), R\})_{i=1}^N, \quad x \in X_n.
\end{equation*}
Since the measure $\nu_n := {\Phi_n}_* m_{X_n}$ is supported in the compact set $\left\{z \in \R^N \midd  \|z\|_{\infty} \leq R \right\}$, a sequence $\{\nu_n\}_{n \in \N}$ is tight. Thus the sequence $\{\nu_n\}$ has a subsequence converging weakly to a Borel probability measure $\nu$ on $\R^N$. We define an mm-space $Z$ as
\begin{equation*}
Z : = (\supp{\nu}, \|\cdot\|_{\infty}, \nu).
\end{equation*}
Since $X_n$ dominates $(\supp{\nu_n}, \|\cdot\|_{\infty}, \nu_n)$ for each $n$, we have 
\begin{equation*}
\lim_{n \to \infty}\square(Z, \Py X_n) = 0.
\end{equation*}

For any $i=1, \ldots, N$, we set
\begin{equation*}
A_i := \left\{(z_i)_{i=1}^N \in \R^N \midd z_i = 0 \text{ and } |z_j| \geq \frac{\delta}{2} \text{ for any } j \neq i\right\}.
\end{equation*}

\begin{claim}\label{claim_inc_con}
For any $i =1, \ldots, N$ and for any sufficiently large $n$, we have
\begin{equation*}
\Phi_n(A_{n, i}) \subset A_i.
\end{equation*}
\end{claim}

\begin{proof}
Take any $i \in \{1, \ldots, N\}$ and fix it. It is sufficient to prove that  
\begin{equation*}
\liminf_{n \to \infty} d_{X_n}(A_{n,i}, A_{n, j}) \geq \delta
\end{equation*}
for every $j \neq i$. We put $\alpha := \liminf_{n \to \infty} d_{X_n}(A_{n,i}, A_{n, j})$. There exists a subsequence $\{n_k\} \subset \{n\}$ such that $d_{X_{n_k}}(A_{n_k,i}, A_{n_k, j}) \to \alpha$ as $k \to \infty$. We find two sequences $\{x_m\}_m \subset A_{n_k, i}$ and $\{x'_m\}_m \subset A_{n_k, j}$ such that $d_{X_{n_k}}(x_m, x'_m) \to d_{X_{n_k}}(A_{n_k,i}, A_{n_k, j})$ as $m \to \infty$. Then we have
\begin{align*}
F(\delta) &\leq F(d_Y(B_i, B_j)) \leq F(d_{Y}(f_n(x_m), f_n(x'_m))) \\
&\leq F_{n_k}(d_{X_{n_k}}(x_m, x'_m)) + \varepsilon_{n_k} \\
& \to  F_{n_k}(d_{X_{n_k}}(A_{n_k, i}, A_{n_k, j})) + \varepsilon_{n_k} \text{ as }  m \to \infty
\end{align*}
for each $k$. Moreover, by Lemma \ref{ptcpt}, we have
\begin{equation*}
\lim_{k \to \infty} F_{n_k}(d_{X_{n_k}}(A_{n_k, i}, A_{n_k, j})) = F(\alpha).
\end{equation*}
Combining these implies $F(\delta) \leq F(\alpha)$, so that $\delta \leq \alpha$ since $F$ is increasing. This completes the proof.
\end{proof}

We define a subset $\widetilde{Z} \subset Z$ by
\begin{equation*}
\widetilde{Z} := \supp{\nu} \cap \bigsqcup_{i=1}^N A_i
\end{equation*}
and define a map $g \colon Z \to \dot{Y}$ by
\begin{equation*}
g(z) := \left\{ \begin{array}{ll} y_i & \text{ if } z \in A_i, \\ y_0 & \text{ if } z \in Z \setminus \widetilde{Z}. \end{array}\right.
\end{equation*}
If $B_0 = \emptyset$, then we set $g(z)$, for $z \in Z \setminus \widetilde{Z}$, an arbitrary point of $\dot{Y}$.

\begin{claim}\label{claim_1-Lip_con}
The map $g \colon Z \to \dot{Y}$ is $1$-Lipschitz up to $2\varepsilon$ and fulfills
\begin{equation*}
\prok(g_* \nu, m_{\dot{Y}}) \leq \varepsilon.
\end{equation*}
\end{claim}

\begin{proof}
The proof of $\prok(g_* \nu, m_{\dot{Y}}) \leq \varepsilon$ is completely same as in the proof of Claim \ref{claim_1-Lip}. We prove that the map $g \colon Z \to \dot{Y}$ is $1$-Lipschitz up to $2\varepsilon$. Since
\begin{equation*}
\nu(Z \setminus \tilde{Z}) = g_* \nu(\{y_0\}) \leq m_{\dot{Y}}(\{y_0\}) \leq \varepsilon,
\end{equation*}
we see that $\nu(\tilde{Z}) \geq 1- \varepsilon$. It is sufficient to prove that
\begin{equation}\label{claim_1-Lip_con:eq}
d_{\dot{Y}}(g(z), g(z')) \leq \|z - z'\|_\infty + 2\varepsilon
\end{equation}
holds for any $z, z' \in \tilde{Z}$. We take any $i, j \in \{1, \ldots, N\}$ and fix them. We take any points $z \in A_i$ and $z' \in A_j$. There exist $x_n \in A_{n, i}$ and $x'_n \in A_{n, j}$, $n =1, 2, \ldots$, such that 
\begin{equation*}
\|\Phi_n(x_n) - z \|_\infty, \|\Phi_n(x'_n) - z' \|_\infty \to 0 \text{ as } n \to \infty.
\end{equation*}
Since we have
\begin{equation*}
F(d_Y(y_i, y_j) -2\varepsilon) \leq F(d_Y(B_i, B_j)) \leq F(d_Y(f_n(x_n), f_n(\tilde{x}'))) \leq F_n(d_{X_n}(x_n, \tilde{x}')) + \varepsilon_n
\end{equation*}
for any $\tilde{x}' \in A_{n, j}$, we obtain
\begin{equation*}
F(d_Y(y_i, y_j) -2\varepsilon) \leq \limsup_{n\to\infty} F_n(d_{X_n}(x_n, A_{n, j})).
\end{equation*}
On the other hand, we see that
\begin{equation*}
F(d_Y(y_i, y_j) -2\varepsilon)  \leq F(R) = \lim_{n \to \infty} F_n(R).
\end{equation*}
Combining these implies that
\begin{align*}
F(d_Y(y_i, y_j) -2\varepsilon) &\leq \limsup_{n\to\infty} F_n(\min\{d_{X_n}(x_n, A_{n, j}), R\}) + 3\varepsilon \\
& =  \limsup_{n\to\infty} F_n((\Phi_n(x_n))_j) \\
& = \limsup_{n\to\infty} F_n((\Phi_n(x_n) - \Phi_n(x'_n))_j).
\end{align*}
By \eqref{liminf:eq} and Lemma \ref{liminf:lem}, we have 
\begin{equation*}
F(d_Y(y_i, y_j) -2\varepsilon) \leq \limsup_{n\to\infty} F_n(\|\Phi_n(x_n) - \Phi_n(x'_n)\|_\infty).
\end{equation*}
Furthermore, since  $\|\Phi_n(x_n) - \Phi_n(x'_n)\|_\infty \to \|z - z'\|_\infty$ as $n \to \infty$, we have
\begin{equation*}
F(d_Y(y_i, y_j) - 2\varepsilon)  \leq F(\|z - z'\|_\infty),
\end{equation*}
which leads to
\begin{equation*}
d_{\dot{Y}}(g(z), g(z')) = d_Y(y_i, y_j)  \leq \|z - z'\|_\infty + 2\varepsilon
\end{equation*}
since $F$ is increasing. Therefore we obtain \eqref{claim_1-Lip_con:eq}.
\end{proof}
By Claim \ref{claim_1-Lip_con}, the composition $h := \iota \circ g \colon Z \to Y$ is $1$-Lipschitz up to $3\varepsilon$ and $\prok(h_* m_Z, m_Y) \leq 2\varepsilon$ holds. The proof is completed. 
\end{proof}

\begin{lem}\label{lem_inj}
Let $F \colon [0,+\infty) \to [0,+\infty)$ be a continuous nondecreasing metric preserving function. Then, the following {\rm (1) -- (3)} are equivalent to each other.
\begin{enumerate}
\item $F$ is increasing.
\item If two pyramids $\Py$, $\PyQ \in \Pi$ satisfy $F(\PyQ) \subset F(\Py)$, then $\PyQ \subset \Py$.
\item The map $\widetilde{F} \colon \Pi \ni \Py \mapsto F(\Py) \in \Pi$ is injective.
\end{enumerate}
\end{lem}

\begin{proof}
`(2) $\Rightarrow$ (3)' is obvious. 

We first prove `(3) $\Rightarrow$ (1)'. Suppose that $F$ is not increasing. By $F(0) = 0$ and the continuity of $F$, there exist two positive real numbers $s < s'$ such that $F(s) = F(s')$. Two mm-spaces $X$ and $Y$ are defined as 
\begin{equation*}
X := (\{0, s\}, |\cdot|, \delta_{0,s}) \quad \text{ and } \quad Y := (\{0, s'\}, |\cdot|, \delta_{0,s'}),
\end{equation*}
where the notation is same as in the proof of Proposition \ref{Py_ptwise}. Then, $X$ and $Y$ are not mm-isomorphic to each other, but $F(X) = F(Y)$, so that $\widetilde{F}$ is not injective.

We next prove `(1) $\Rightarrow$ (2)'. Assume that two pyramids $\Py$, $\PyQ \in \Pi$ satisfy $F(\PyQ) \subset  F(\Py)$. We take any mm-space $Y \in \PyQ$. Since $F(Y) \in F(\PyQ) \subset F(\Py)$, there exist mm-spaces $X_n \in \Py$, $n=1, 2, \ldots$, such that 
\begin{equation*}
\lim_{n \to \infty}\square(F(Y), \Py{F(X_n)}) = 0.
\end{equation*}
By Lemma \ref{lip_up_con}, we obtain
\begin{equation*}
\square(Y, \Py) \leq \lim_{n \to \infty} \square(Y, \Py X_n) = 0,
\end{equation*}
which implies $Y \in \Py$. This completes the proof of this lemma.
\end{proof}

\begin{prop}\label{Fapprox}
Let $F_n, F \colon [0, +\infty) \to [0, +\infty)$ be continuous metric preserving functions. Assume that $F_n$ converges pointwise to $F$ as $n \to \infty$ and that $F$ is nondecreasing.  Then, for any pyramid $\Py \in \Pi$, there exists a sequence $\{Y_n\}_{n \in \N}$ of mm-spaces such that $Y_n$ converges weakly to $\Py$ and $F_n(Y_n)$ converges weakly to $F(\Py)$.
\end{prop}

\begin{proof}
We take any pyramid $\Py \in \Pi$. Let $\{Y_m\}_{m \in \N}$ be an approximation of $\Py$ (see Lemma \ref{approximation}). Note that $Y_m$ converges weakly to $\Py$ as $m \to \infty$. For each $m$, by \cite{prod}*{Corollary 4.4} (see Remark \ref{sum_table}), we have 
\begin{equation*}
\lim_{n \to \infty} \square(F_n(Y_m), F(Y_m)) = 0.
\end{equation*}
Choosing sufficiently small $m$ than $n$, there exists a sequence $\{m(n)\}_{n \in \N}$ with $m(n) \to \infty$ as $n \to \infty$ such that
\begin{equation}\label{approx:eq}
\lim_{n \to \infty} \square(F_n(Y_{m(n)}), F(Y_{m(n)})) = 0.
\end{equation}
We prove that $F_n(Y_{m(n)})$ converges weakly to $F(\Py)$ as $n \to \infty$.

Take any mm-space $Z \in F(\Py)$ and any real number $\varepsilon > 0$. There exist mm-spaces $X \in \Py$ and $Z' \in \X$ such that 
\begin{equation}\label{approx:liminfeq1}
Z' \prec F(X) \quad \text{ and } \quad \square(Z, Z') < \varepsilon.
\end{equation}
Moreover, there exist mm-spaces $X_n$, $n = 1, 2, \ldots$, such that $Y_{m(n)}$ dominates $X_n$ for all $n$ and $\square(X_n, X) \to 0$ as $n \to \infty$. Since $F$ is nondecreasing, we see that 
\begin{equation}\label{approx:liminfeq2}
F(X_n) \prec F(Y_{m(n)}) \text{ for every } n \quad \text{ and } \quad \lim_{n \to \infty} \square(F(X_n), F(X)) = 0.
\end{equation}
Combining \eqref{approx:eq}, \eqref{approx:liminfeq1}, \eqref{approx:liminfeq2}, and Lemma \ref{lem4.22} means that, for any sufficiently large $n$, there exist maps $f_n \colon F_n(Y_{m(n)}) \to Z$ such that $f_n$ is $1$-Lipschitz up to $4\varepsilon$ and $\prok({f_n}_* m_{Y_{m(n)}}, m_Z) \leq 4\varepsilon$. Thus, by Lemma \ref{Py:lem}, we have
\begin{equation*}
\square(Z, \Py{F_n(Y_{m(n)})}) \leq 16\varepsilon
\end{equation*}
for any sufficiently large $n$. We obtain $\lim_{n \to \infty} \square(Z, \Py{F_n(Y_{m(n)})}) = 0$.

On the other hand, we take any mm-space $Z \in \X$ such that 
\begin{equation*}
\liminf_{n \to \infty} \square(Z, \Py{F_n(Y_{m(n)})}) = 0.
\end{equation*}
Taking a subsequence of $n$, we can assume that $\square(Z, \Py{F_n(Y_{m(n)})}) \to 0$ as $n \to \infty$. By Lemma \ref{Py:lem}, \eqref{approx:eq}, and Lemma \ref{lem4.22}, there exist maps $g_n \colon F(Y_{m(n)}) \to Z$ and a sequence $\varepsilon_n \to 0$ such that $g_n$ is $1$-Lipschitz up to $\varepsilon_n$ and $\prok({g_n}_* m_{Y_{m(n)}}, m_Z) \leq \varepsilon_n$ holds for every $n$. Since $F(Y_{m(n)}) \in F(\Py)$ follows from $Y_{m(n)} \in \Py$, we have $Z \in F(\Py)$ by Corollary \ref{Py:lem_cor}.

Therefore $F_n(Y_{m(n)})$ converges weakly to $F(\Py)$ as $n \to \infty$. The proof is completed.
\end{proof}

\begin{proof}[Proof of {\rm `}$(1) \Rightarrow (2)${\rm '} of Theorem \ref{Py_con}]
Assume the condition (1). We first prove \eqref{ptwise}. By Proposition \ref{ptwise_con}, it is sufficient to prove that if $F_n(s_n) \to F(s)$ for given positive real numbers $s_n$, $n = 1, 2, \ldots$, and $s$, then we have $s_n \to s$. Let $s_n$, $n = 1, 2, \ldots$, and $s$ be positive real numbers with $F_n(s_n) \to F(s)$. We define mm-spaces $X_n$ and $X$ as
\begin{equation*}
X_n := (\{0, s_n\}, |\cdot|, \delta_{0,s_n}) \quad \text{and} \quad X := (\{0, s\}, |\cdot|, \delta_{0,s}),
\end{equation*}
where the notation is same as in the proof of Proposition \ref{Py_ptwise}. Note that $F_n(X_n)$ $\square$-converges to $F(X)$. By the condition (1), $X_n$ concentrates to $X$. Thus we have $s_n \to s$ as $n \to \infty$ in the same way as the proof of Proposition \ref{Py_ptwise}. We obtain \eqref{ptwise}. Moreover, \eqref{liminf:eq} follows from Lemma \ref{prod:claim5.1} directly.

We next prove \eqref{incr}. By Lemma \ref{lem_inj}, it is sufficient to prove that the map $\Py \mapsto F(\Py)$ is injective. We take any two pyramids $\Py$, $\PyQ$ and assume that $F(\Py) = F(\PyQ)$. By Proposition \ref{Fapprox}, there exists a sequence $\{Y_n\}_{n \in \N}$ of mm-spaces such that $Y_n$ converges weakly to $\Py$ and $F_n(Y_n)$ converges weakly to $F(\Py)$. It follows from $F(\Py) = F(\PyQ)$ and the condition (1) that $Y_n$ converges weakly to $\PyQ$, which implies that $\Py = \PyQ$. This completes the proof.
\end{proof}

We next prove {\rm`(2) $\Rightarrow$ (1)'} of Theorem \ref{Py_con}. Let $F_n$, $F$ be continuous metric preserving functions satisfying \eqref{ptwise} \eqref{liminf:eq}, and \eqref{incr}. Take any sequence $\{X_n\}_{n \in \N}$ of mm-spaces and any pyramid $\Py \in \Pi$ such that $F_n(X_n)$ converges weakly to $F(\Py)$. Our goal is to prove that $X_n$ converges weakly to $\Py$ as $n \to \infty$.

\begin{prop}\label{liminf:prop_con}
For any $Y \in \Py$, we have 
\begin{equation*}
\lim_{n \to \infty}\square(Y, \Py{X_n}) = 0.
\end{equation*}
\end{prop}

\begin{proof}
We take any mm-space $Y \in \Py$. Since $F(Y) \in F(\Py)$ and $F_n(X_n)$ converges weakly to $F(\Py)$, we see that
\begin{equation*}
\lim_{n \to \infty}\square(F(Y), \Py{F_n(X_n)}) = 0.
\end{equation*}
By Lemma \ref{lip_up_con}, we have
\begin{equation*}
\lim_{n \to \infty}\square(Y, \Py{X_n}) = 0.
\end{equation*}
The proof is completed.
\end{proof}

\begin{prop}\label{limsup:prop_con}
If an mm-space $Y$ satisfies
\begin{equation*}
\liminf_{n \to \infty} \square(Y, \Py X_n) = 0,
\end{equation*}
then $Y \in \Py$.
\end{prop}

\begin{proof}
Choosing a subsequence of $n$, we can assume that $\square(Y, \Py{X_n}) \to 0$ as $n \to \infty$. Then, by Lemma \ref{Py:lem}, there exist Borel measurable maps $f_n \colon X_n \to Y$ and a sequence $\varepsilon_n \to 0$ such that $f_n$ is $1$-Lipschitz up to $\varepsilon_n$ and $\prok^{d_Y}({f_n}_* m_{X_n}, m_Y) \leq \varepsilon_n$ holds for every $n$. By Lemma \ref{lip_up}, there exists a sequence $\delta_n \to 0$ such that $f_n \colon F_n(X_n) \to F(Y)$ is $1$-Lipschitz up to $\delta_n$ and $\prok^{F\circ d_Y}({f_n}_* m_{X_n}, m_Y) \leq \delta_n$ holds. By Lemma \ref{Py:lem}, we have
\begin{equation*}
\square(F(Y), \Py{F_n(X_n)}) \leq 4\delta_n \to 0 \text{ as } n \to \infty
\end{equation*}
Since $F_n(X_n)$ converges weakly to $F(\Py)$, we have $F(Y) \in F(\Py)$. Thus, since $F$ is increasing, we obtain $Y \in \Py$. The proof is completed.
\end{proof}

\begin{proof}[Proof of {\rm `}$(2) \Rightarrow (1)${\rm '} of Theorem \ref{Py_con}]
This follows from Proposition \ref{liminf:prop_con} and Proposition \ref{limsup:prop_con}.
\end{proof}

\begin{proof}[Proof of Corollary \ref{Py:cor} {\rm (B)}]
`(B1) $\Rightarrow$ (B2)' follows from Theorem \ref{Py_con} directly. We prove `(B2) $\Rightarrow$ (B1)'. Assume that $F_n(\Py_n)$ converges weakly to $F(\Py)$.

We take any mm-space $Y \in \Py$. Since $F(Y) \in F(\Py)$, we have
\begin{equation*}
\lim_{n \to \infty} \square(F(Y), F_n(\Py_n)) = 0.
\end{equation*}
Then, there exist mm-spaces $X_n \in \Py_n$, $n = 1, 2, \ldots$, such that
\begin{equation*}
\lim_{n \to \infty} \square(F(Y), \Py F_n(X_n)) = 0.
\end{equation*}
By Lemma \ref{lip_up_con}, we have
\begin{equation*}
\limsup_{n \to \infty} \square(Y, \Py_n) \leq \lim_{n \to \infty} \square(Y, \Py X_n) = 0.
\end{equation*}

On the other hand, we assume that an mm-space $Y$ satisfies
\begin{equation*}
\liminf_{n \to \infty} \square(Y, \Py_n) = 0.
\end{equation*}
Taking a subsequence of $n$, we can assume that there exist mm-spaces $X_n \in \Py_n$ and $Y_n \in \X$, $n = 1, 2, \ldots$, such that $X_n$ dominates $Y_n$ for every $n$ and $\square(Y_n, Y) \to 0$ as $n \to \infty$. Since all $F_n$ are nondecreasing, we see that $F_n(X_n)$ dominates $F_n(Y_n)$ for every $n$ and $\square(F_n(Y_n), F(Y)) \to 0$ as $n \to \infty$. Thus, we have
\begin{equation*}
\limsup_{n \to \infty} \square(F(Y), F_n(\Py_n)) \leq \lim_{n \to \infty} \square(F(Y), F_n(Y_n)) = 0,
\end{equation*}
which implies $F(Y) \in F(\Py)$, so that $Y \in\Py$ since $F$ is increasing.

Combining these means that $\Py_n$ converges weakly to $\Py$. This completes the proof of this corollary.
\end{proof}

\section{Box-convergence/concentration of metric transformed spaces}

The goal in this section is to prove Theorem \ref{mm_con}.

\begin{prop}\label{liminf_ness}
Let $F_n, F \colon [0, +\infty) \to [0, +\infty)$, $n = 1, 2, \ldots$, be continuous metric preserving functions. Assume that $F_n$ converges pointwise to $F$ as $n \to \infty$. If there exists a sequence $s_n \to \infty$ such that
\begin{equation*}
\liminf_{n \to \infty} F_n(s_n) < \sup{F},
\end{equation*}
then there exist a sequence $\{X_n\}_{n \in \N}$ of mm-spaces such that $F_n(X_n)$ $\square$-converges to $F(X)$ for an mm-space $X$ but $X_n$ does not concentrate.
\end{prop}

\begin{proof}
We take a sequence $s_n \to \infty$ such that
\begin{equation*}
\liminf_{n \to \infty} F_n(s_n) < \sup{F}.
\end{equation*}
We set $\alpha := \liminf_{n \to \infty} F_n(s_n)$. There exists a real number $\beta$ such that $F(\beta) = \alpha$. We find a subsequence $\{n_i\} \subset \{n\}$ such that
\begin{equation*}
\lim_{i \to \infty} F_{n_i}(s_{n_i}) = \alpha = F(\beta).
\end{equation*}
We define an mm-space
\begin{equation*}
X := (\{0, \beta\}, |\cdot|, \delta_{0, \beta}),
\end{equation*}
and define mm-spaces
\begin{equation*}
X_n := \left\{\begin{array}{ll} (\{0, s_{n_i}\}, |\cdot|, \delta_{0, s_{n_i}}) & \text{ if } n = n_i, \\ X & \text{ otherwise,} \end{array}\right.
\end{equation*}
where the notation is same as in the proof of Proposition \ref{Py_ptwise}. It is easy to see that $F_n(X_n)$ $\square$-converges to $F(X)$. Moreover, since $s_{n_i} \to \infty$, we see that $X_{n_i}$ converges weakly to the pyramid
\begin{equation*}
\Py := \left\{ (\{0, s\}, |\cdot|, \delta_{0, s}) \midd s \geq 0 \right\}
\end{equation*}
as $i \to \infty$. Since $\Py$ is not $\square$-precompact, $X_{n_i}$ does not concentrate to any mm-spaces. This completes the proof.
\end{proof}

The following proposition is a corollary of Lemma \ref{lem_inj}.

\begin{prop}\label{prop_inj}
Let $F \colon [0,+\infty) \to [0,+\infty)$ be a continuous metric preserving function {\rm (}we do not assume that $F$ is nondecreasing{\rm )}. Then, the following {\rm (1)} and {\rm (2)} are equivalent to each other.
\begin{enumerate}
\item $F$ is increasing.
\item The map $\X \ni X \mapsto F(X) \in \X$ is injective.
\end{enumerate}
\end{prop}

\begin{proof}
The implication `(1) $\Rightarrow$ (2)' follows from Lemma \ref{lem_inj} directly. The proof of `(2) $\Rightarrow$ (1)' is completely same as the proof of `(3) $\Rightarrow$ (1)' of Lemma \ref{lem_inj}. Note that we do not need the assumption that $F$ is nondecreasing in this proof. 
\end{proof}

\begin{proof}[Proof of Theorem \ref{mm_con}]
We first prove `(2) $\Rightarrow$ (1)' and `(3) $\Rightarrow$ (1)' together. Assume that one of the conditions (2) or (3) holds. We have \eqref{ptwise} by the completely same discussion as the proof of Theorem \ref{Py_con}. Moreover, \eqref{liminf:eq} follows from Proposition \ref{liminf_ness} and Lemma \ref{liminf:lem}. We verify \eqref{incr}. We take any two mm-spaces $X$ and $Y$ with $F(X) = F(Y)$. By \cite{prod}*{Corollary 4.4} (see also Remark \ref{sum_table}), we see that $F_n(X)$ $\square$-converges to $F(X)$. By $F(X) = F(Y)$ and our assumption, we have $X = Y$. Thus the map $\X \ni X \mapsto F(X) \in \X$ is injective. Combining this and Proposition \ref{prop_inj} implies that $F$ is increasing.  Thus we obtain `(2) $\Rightarrow$ (1)' and `(3) $\Rightarrow$ (1)'.

The implication `(1) $\Rightarrow$ (3)' follows from Theorem \ref{Py_con} directly. We next prove  `(1) $\Rightarrow$ (2)'. Assume that \eqref{ptwise}, \eqref{liminf:eq} and \eqref{incr} hold. Take any sequence $\{X_n\}_{n \in \N}$ of mm-spaces and any mm-space $X$ such that $F_n(X_n)$ $\square$-converges to $F(X)$. We prove that $X_n$ $\square$-converges to $X$. 

Take any real number $\varepsilon > 0$. By the inner regularity of $m_X$, there exists a compact set $K \subset X$ such that $m_X(K) \geq 1 - \varepsilon$. Let
\begin{equation*}
D_\varepsilon := \left\{\begin{array}{ll} F(\diam{K}) + 4\varepsilon & \text{ if } \sup{F} = +\infty, \\ \frac{1}{2}(F(\diam{K}) + \sup{F}) & \text{ if } \sup{F} < +\infty. \end{array}\right.
\end{equation*}
Note that $D_\varepsilon < +\infty$ and $[0, D_\varepsilon] \subset \mathrm{Im}F$. Since $F$ is continuous and increasing, so is the inverse function $F^{-1} \colon \mathrm{Im}F \to [0, +\infty)$. Moreover, $F^{-1}$ is uniformly continuous on $[0, D_\varepsilon]$. Let $\omega_\varepsilon$ is the minimal modulus of continuity of $F^{-1}|_{[0, D_\varepsilon]}$, that is,
\begin{equation*}
\omega_\varepsilon(\delta) := \sup\left\{ |F^{-1}(s) - F^{-1}(t)| \midd s, t \in [0, D_\varepsilon] \text{ with } |s - t| \leq \delta\right\}.
\end{equation*}
We take any real number $\delta$ such that
\begin{equation*}
0 < \delta < \min\left\{\varepsilon, \frac{\sup{F} - F(\diam{K})}{8}\right\} \quad \text{ and } \quad \omega_\varepsilon(\delta) < \varepsilon.
\end{equation*}
Since $F_n(X_n)$ $\square$-converges to $F(X)$, there exist $\delta$-mm-isomorphisms $f_n \colon F_n(X_n) \to F(X)$ for sufficiently large $n$.

\begin{claim}\label{mm_con:claim}
For every sufficiently large $n$, the map $f_n \colon X_n \to X$ is a $3\varepsilon$-mm-isomorphism.
\end{claim}

\begin{proof}[Proof of Claim \ref{mm_con:claim}]
We first prove $\prok^{d_X}({f_n}_* m_{X_n}, m_X) \leq \varepsilon$. For any subset $A \subset X$, we have
\begin{equation*}
U_\delta^{F \circ d_X}(A) \subset U_{F^{-1}(\delta)}^{d_X}(A).
\end{equation*}
In fact, taking any point $y \in U_\delta^{F \circ d_X}(A)$, it holds that $F(d_X(y, A)) < \delta$, which implies that $d_X(y, A) < F^{-1}(\delta)$. Combining $\prok^{F \circ d_X}({f_n}_* m_{X_n}, m_X) \leq \delta$ and this leads to
\begin{equation*}
m_X(A) \leq {f_n}_* m_{X_n}(U_\delta^{F \circ d_X}(A)) + \delta \leq {f_n}_* m_{X_n}(U_{F^{-1}(\delta)}^{d_X}(A)) + \delta,
\end{equation*}
which implies that $\prok^{d_X}({f_n}_* m_{X_n}, m_X) \leq \max\{F^{-1}(\delta), \delta\}$. Since 
\begin{equation*}
F^{-1}(\delta) = |F^{-1}(\delta) - F^{-1}(0)| \leq \omega_\varepsilon(\delta) < \varepsilon,
\end{equation*}
we obtain $\prok^{d_X}({f_n}_* m_{X_n}, m_X) \leq \varepsilon$.

Let $X'_n \subset F_n(X_n)$ be a nonexceptional domain of $f_n \colon F_n(X_n) \to F(X)$ and let
\begin{equation*}
\widetilde{X}_n := X'_n \cap f_n^{-1}(U_\delta^{F \circ d_X}(K)).
\end{equation*}
We see that
\begin{equation*}
m_{X_n}(\widetilde{X}_n) \geq m_{X_n}(X'_n) + {f_n}_* m_{X_n}(U_\delta^{F \circ d_X}(K)) - 1 \geq m_{X}(K) -2\delta \geq 1 - 3\varepsilon.
\end{equation*}
It is sufficient to prove that we have
\begin{equation*}
|d_{X_n}(x, x') - d_X(f_n(x), f_n(x'))| \leq 2\varepsilon
\end{equation*}
for every sufficiently large $n$ and for any $x, x' \in \widetilde{X}_n$. We see that 
\begin{equation*}
F(d_X(f_n(x), f_n(x'))) \leq D_\varepsilon, \quad F_n(d_{X_n}(x, x')) \leq D_\varepsilon, \quad \text{and} \quad F(d_{X_n}(x, x')) \leq D_\varepsilon
\end{equation*}
for every sufficiently large $n$ and for any $x, x' \in \widetilde{X}_n$. In fact, we have
\begin{equation*}
F_n(d_{X_n}(x, x')) \leq F(d_X(f_n(x), f_n(x'))) + \delta \leq F(\diam{K}) + 3\delta \leq D_\varepsilon.
\end{equation*}
Suppose that $d_{X_n}(x, x') \geq F^{-1}(D_\varepsilon + \eta)$ for some $\eta > 0$. By Lemma  \ref{liminf:lem}, we have
\begin{equation*}
F_n(F^{-1}(D_\varepsilon + \eta)) \leq F_n(d_{X_n}(x, x')) + \delta \leq D_\varepsilon.
\end{equation*}
As $n \to \infty$, this implies the contradiction $D_\varepsilon + \eta \leq D_\varepsilon$. Thus we obtain $d_{X_n}(x, x') \leq F^{-1}(D_\varepsilon)$, that is, $F(d_{X_n}(x, x')) \leq D_\varepsilon$. Combining this and Lemma \ref{ptcpt} implies that, for any $x, x' \in \widetilde{X}_n$,
\begin{align*}
&|d_{X_n}(x, x') - d_X(f_n(x), f_n(x'))| \\
&\leq |d_{X_n}(x, x') - F^{-1}(F_n(d_{X_n}(x, x')))| + |F^{-1}(F_n(d_{X_n}(x, x'))) - d_X(f_n(x), f_n(x'))| \\
&\leq \omega_\varepsilon(|F(d_{X_n}(x, x')) - F_n(d_{X_n}(x, x'))|) + \omega_\varepsilon(|F_n(d_{X_n}(x, x')) - F(d_X(f_n(x), f_n(x')))|) \\
&\leq 2\omega_\varepsilon(\delta) < 2\varepsilon,
\end{align*}
where the third inequality follows from, for every sufficiently large $n$,
\begin{equation*}
\sup_{s \in [0, F^{-1}(D_\varepsilon)]}|F(s) - F_n(s)| \leq \delta.
\end{equation*}
Thus the map $f_n \colon X_n \to X$ is a $3\varepsilon$-mm-isomorphism. The proof of this claim is completed.
\end{proof}
By Claim \ref{mm_con:claim}, we see that $X_n$ $\square$-converges to $X$ as $n \to \infty$. Thus we obtain `(1) $\Rightarrow$ (2)'. This completes the proof of Theorem \ref{mm_con}.
\end{proof}

\section{Application: spheres and projective spaces}

\subsection{Gaussian space}

Let $\lambda$ be a positive real number. The product
\begin{equation*}
\gamma^n_{\lambda^2} := \bigotimes_{i=1}^n \gamma^1_{\lambda^2}
\end{equation*}
of the one-dimensional centered Gaussian measure $\gamma^1_{\lambda^2}$ of variance $\lambda^2$ is the $n$-dimensional centered Gaussian measure on $\R^n$ of variance $\lambda^2$. We call the mm-space 
\begin{equation*}
\Gamma^n_{\lambda^2} := (\R^n,\|\cdot\|,\gamma^n_{\lambda^2})
\end{equation*}
the \emph{$n$-dimensional Gaussian space with variance $\lambda^2$}. 

For $1 \le k \le n$, we denote by $\pi^n_k : \R^n \to \R^k$ the natural projection, that is,
\begin{equation*}
\pi^n_k(x_1,x_2,\dots,x_n) := (x_1,x_2,\dots,x_k), \quad (x_1,x_2,\dots,x_n) \in \R^n.
\end{equation*}
Since the projection $\pi^n_{n-1} : \Gamma^n_{\lambda^2} \to \Gamma^{n-1}_{\lambda^2}$ is $1$-Lipschitz continuous and measure-preserving for any $n \geq 2$, the Gaussian space $\Gamma^n_{\lambda^2}$ is monotone nondecreasing in $n$ with respect to the Lipschitz order, so that, as $n \to \infty$, the Gaussian space $\Gamma^n_{\lambda^2}$ converges weakly to the pyramid
\begin{equation*}
\Py {\Gamma^\infty_{\lambda^2}} := \overline{\bigcup_{n=1}^\infty \Py {\Gamma^n_{\lambda^2}}}^{\, \square}. 
\end{equation*}
We call $\Py {\Gamma^\infty_{\lambda^2}}$ the \emph{virtual Gaussian space with variance $\lambda^2$}.
We remark that the infinite product measure
\begin{equation*}
\gamma^\infty_{\lambda^2} := \bigotimes_{i=1}^\infty \gamma^1_{\lambda^2}
\end{equation*}
is a Borel probability measure on $\R^\infty$ with respect to the product topology, but is not Borel with respect to the $l^2$-norm.

Let $F = \R$, $\C$, or $\Hb$, where $\Hb$ is the algebra of quaternions, and let $d := \dim_\R F$. We consider the Hopf action on $\Gamma^{dn}_{\lambda^2}$ by identifying $\R^{dn}$ with $F^n$. Recall that the Hopf action is the following $U^F(1)$-action on $F^n$: 
\begin{equation*}
F^n \times U^F(1) \ni (z, t) \mapsto zt \in F^n,
\end{equation*}
where $U^F(1) := \left\{t \in F \midd \|t\| = 1 \right\}$ is a group under multiplication. Since the projection $\pi^{dn}_{dk} \colon F^n \to F^k$, $k \leq n$, is $U^F(1)$-equivariant (i.e., 
\begin{equation*}
\pi^{dn}_{dk}(zt) =  \pi^{dn}_{dk}(z)t
\end{equation*}
for any $t \in U^F(1)$ and for any $z \in F^n$), there exists a unique map $\bar{\pi}^{dn}_{dk} \colon {F^n}/{U^F(1)} \to {F^k}/{U^F(1)}$ such that $q \circ \pi^{dn}_{dk} = \bar{\pi}^{dn}_{dk} \circ q$, where $q$ is the quotient map of the Hopf action. The Hopf action is isometric with respect to the Euclidean distance and also preserves the Gaussian measure $\gamma^{dn}_{\lambda^2}$. Let
\begin{equation*}
{\Gamma^{dn}_{\lambda^2}}/{U^F(1)} := ({F^n}/{U^F(1)}, d_{{F^n}/{U^F(1)}}, \bar{\gamma}^{dn}_{\lambda^2})
\end{equation*}
be the quotient space with the induced mm-structure, that is,
\begin{align*}
d_{{F^n}/{U^F(1)}}([z], [w]) &:= \inf_{z' \in [z], w' \in [w]} \|z' - w'\|, \quad [z], [w] \in F^n/U^F(1), \\
\bar{\gamma}^{dn}_{\lambda^2}(A) &:= \gamma^{dn}_{\lambda^2}(\left\{z \in F^n \midd [z] \in A \right\}), \quad A \subset F^n/U^F(1).
\end{align*}
Since the map $\bar{\pi}^{dn}_{d(n - 1)} \colon {\Gamma^{dn}_{\lambda^2}}/{U^F(1)} \to {\Gamma^{d(n - 1)}_{\lambda^2}}/{U^F(1)}$ is 1-Lipschitz continuous and measure-preserving, the quotient space ${\Gamma^{dn}_{\lambda^2}}/{U^F(1)}$ is monotone increasing in $n$ with respect to the Lipschitz order. The Hopf quotient space ${\Gamma^{dn}_{\lambda^2}}/{U^F(1)}$ converges weakly to the pyramid
\begin{equation*}
\Py {\Gamma^\infty_{\lambda^2}}/{U^F(1)} := \overline{\bigcup_{n=1}^\infty \Py {\Gamma^{dn}_{\lambda^2}/{U^F(1)}}}^{\, \square}. 
\end{equation*}

\subsection{Weak convergence of spheres and projective spaces}
Let $S^n(r)$ be the $n$-dimensional sphere in $\R^{n+1}$ centered at the origin and of radius $r > 0$. We equip $S^n(r)$ with the standard Riemannian distance function $d_{S^n(r)}$ or the restriction of the Euclidean distance function $\|\cdot\|$. We also equip $S^n(r)$ with the Riemannian volume measure $\sigma^n$ normalized as $\sigma^n(S^n(r)) = 1$. Then $S^n(r)$ is an mm-space. We consider the Hopf quotient
\begin{equation*}
FP^n(r) := {S^{d(n+1)-1}(r)}/{U^F(1)}
\end{equation*}
that has a natural mm-structure induced from that of $S^{d(n+1)-1}(r)$ by the same way as above. This is topologically an $n$-dimensional projective space over $F$. Note that, if $F = \C$ and if the distance function on $S^{2n+1}(r)$ is assumed to be Riemannian, then the distance function on ${\C}P^n(r)$ coincides with that induced from the Fubini-Study metric scaled with factor $r$.

\begin{thm}[\cite{MML}*{Theorem 8.1.1}, \cite{ST}*{Corollary 1.3}]
Let $\{r_n\}_{n=1}^\infty$ be a given sequence of positive real numbers, and let $\lambda_n := r_n /\sqrt{n}$ {\rm (}resp.~$\lambda_n := r_n /\sqrt{dn}${\rm )}. Then we have the following {\rm (1)} and {\rm (2)}.
\begin{enumerate}
\item $\{S^n(r_n)\}_{n \in \N}$ {\rm (}resp.~$\{FP^n(r_n)\}_{n \in \N}${\rm )} is L\'evy family {\rm (}i.e., concentrating to a one point space{\rm )} if and only if $\lambda_n$ converges to $0$ as $n \to \infty$.
\item  $\{S^n(r_n)\}_{n \in \N}$ {\rm (}resp.~$\{FP^n(r_n)\}_{n \in \N}${\rm )} infinitely dissipates if and only if $\lambda_n$ diverges to infinity as $n \to \infty$.
\end{enumerate}
\end{thm}

We omit to state the definition of the infinite dissipation. Dissipation is the opposite notion to concentration. The above theorem claims that the critical scale order for concentration is $\sqrt{n}$. Moreover, in the Euclidean case, the limit of spheres and projective spaces with the critical scale order is known. 

\begin{thm}[\cite{MML}*{Theorem 8.1.1}, \cite{ST}*{Theorem 1.2}] \label{Sthm}
Let $\{r_n\}_{n=1}^\infty$ be a given sequence of positive real numbers, and let $\lambda_n := r_n /\sqrt{n}$ {\rm (}resp. $\lambda_n := r_n /\sqrt{dn}${\rm )}. Assume that $S^n(r_n)$ and $FP^n(r_n)$ have the Euclidean distance function. As $n \to \infty$, $S^n(r_n)$ {\rm (}resp.~$FP^n(r_n)${\rm )} converges weakly to $\Py{\Gamma^\infty_{\lambda^2}}$ {\rm (}resp.~$\Py{\Gamma^\infty_{\lambda^2}}/U^F(1)${\rm )} if and only if $\lambda_n$ converges to a positive real number $\lambda$. 
\end{thm}

\begin{rem}
The `only if' part of the above theorem can be easily checked as follows. In \cite{OS}, the $\kappa$-observable diameter of a pyramid, which is a fundamental invariant of a pyramid, is introduced and the limit formula is proved. The $\kappa$-observable diameter of $\Py{\Gamma^\infty_{\lambda^2}}$ {\rm (}resp.~$\Py{\Gamma^\infty_{\lambda^2}}/U^F(1)${\rm )} is proportional to $\lambda$, so that the pyramid $\Py{\Gamma^\infty_{\lambda^2}}$ {\rm (}resp. $\Py{\Gamma^\infty_{\lambda^2}}/U^F(1)${\rm )} is different for each $\lambda$.
\end{rem}

Roughly speaking, this theorem has been obtained by proving
\begin{itemize}
\item the limit of $S^n(r_n)$ (resp.~$FP^n(r_n)$) dominates $\Py{\Gamma^\infty_{\lambda^2}}$ (resp.~$\Py{\Gamma^\infty_{\lambda^2}}/U^F(1)$),
\item $\Py{\Gamma^\infty_{\lambda^2}}$ (resp.~$\Py{\Gamma^\infty_{\lambda^2}}/U^F(1)$) dominates the limit of $S^n(r_n)$ (resp.~$FP^n(r_n)$)
\end{itemize}
from the construction of maps for the domination. We remark that it would be difficult to find a map for the domination directly in the Riemannian case. Our goal is to prove Theorem \ref{SCP} using the convergence of metric transformed spaces.

\begin{proof}[Proof of Theorem \ref{SCP}]
Let $\{r_n\}_{n=1}^\infty$ be a given sequence of positive real numbers. We define metric preserving functions $F_n \colon  [0, +\infty) \to  [0, +\infty)$, $n = 1, 2, \ldots$, by
\begin{equation*}
F_n(s) := \left\{\begin{array}{ll} 2r_n \sin{\frac{s}{2r_n}} & \text{if } s \leq \pi r_n, \\ 2r_n & \text{if } s > \pi r_n. \end{array}\right.
\end{equation*}
We see that 
\begin{equation*}
\| x - x' \|  = F_n(d_{S^n(r_n)}(x, x'))
\end{equation*}
for any $x, x' \in S^n(r_n)$. Similarly, we have
\begin{equation*}
\inf_{z' \in [z], w' \in [w]}\| z' - w' \|  = F_n(\inf_{z' \in [z], w' \in [w]} d_{S^{d(n+1)-1}(r_n)}(z', w'))
\end{equation*}
for any $[z], [w] \in FP^n(r_n)$. Since $F_n$ converges pointwise to the identity function $s \mapsto s$ as $n \to \infty$, combining Theorem \ref{Sthm} and Theorem \ref{Py_con} implies that $S^n(r_n)$ converges weakly to $\Py{\Gamma^\infty_{\lambda^2}}$ and that $FP^n(r_n)$ converges weakly to $\Py{\Gamma^\infty_{\lambda^2}}/U^F(1)$ as $n \to \infty$ for the Riemannian distance functions.
\end{proof}

\end{document}